\documentclass[preprint,12pt]{elsarticle}

\usepackage{bm}
\usepackage{hyperref}
\usepackage{graphicx}
\usepackage{subfigure}
\usepackage{caption}
\usepackage{overpic}
\usepackage{epsfig}
\usepackage{xcolor}
\usepackage{indentfirst}
\usepackage{fancyhdr}
\usepackage{comment}
\usepackage{makecell, rotating}
\usepackage{float}
\usepackage{fullpage} 
\usepackage{multirow}
\usepackage{booktabs}
\usepackage{xkeyval}
\usepackage{mathtools}
\usepackage{exscale}
\usepackage{relsize}

\usepackage[ruled,linesnumbered]{algorithm2e}
\newcommand{\bbR}{\mathbb{R}}

\newcommand{\bbS}{\mathcal{S}}
\newcommand{\br}{\bm{r}}
\newcommand{\bu}{\mathbf{u}}
\newcommand{\bx}{\bm{x}}

\newcommand{\bc}{\bm{c}}

\newcommand{\ba}{\bm{a}}
\newcommand{\bb}{\bm{b}}

\newcommand{\bv}{\bm{v}}

\newcommand{\bk}{\mathbf{k}}

\newcommand{\calB}{\mathcal{B}}

\newcommand{\calW}{\mathcal{W}}

\newcommand{\bS}{\boldsymbol{S}}
\newcommand{\bM}{\boldsymbol{M}}
\newcommand{\bI}{\boldsymbol{I}}

\DeclareMathOperator*{\argmin}{\mathrm{argmin}}

\usepackage{amssymb}
\usepackage{amsmath}

\journal{*}
 
\begin{document}

\begin{frontmatter}



\title{High-accurate and efficient numerical algorithms for the self-consistent field theory of liquid-crystalline polymers}


\author[1]{Zhijuan He} \author[1]{Kai Jiang\corref{cor}}   \author[1,2]{Liwei Tan}  \author[1]{Xin Wang}

\cortext[cor]{Corresponding author. Email: kaijiang@xtu.edu.cn}

\address[1]{Hunan Key Laboratory for Computation and Simulation in Science and Engineering, Key Laboratory of Intelligent Computing and Information Processing of Ministry of Education, School of Mathematics and Computational Science, Xiangtan University, Xiangtan, Hunan, 411105, China. }
\address[2]{School of Mathematical Sciences, 
	Shanghai Jiao Tong University, 
	Shanghai 200240, P.R. China. }


\begin{abstract}
Self-consistent field theory (SCFT) 
is one of the most widely-used framework in studying the equilibrium phase behaviors of inhomogenous polymers. For liquid crystalline polymeric systems, the main numerical challenges of solving SCFT encompass efficiently solving plenty of six dimensional partial differential equations (PDEs), precisely determining the subtle energy difference among self-assembled structures, and developing effective iterative methods for nonlinear SCF iteration.
To address these challenges, 
this work introduces a suite of high-order and efficient numerical methods tailored for SCFT of liquid-crystalline polymers. These methods include various advaced PDE solvers, an improved Anderson iteration algorithm to accelerate SCFT calculations, and an optimization technique of adjusting the computational domain during the SCF iterations. Extensive numerical tests demonstrate the efficiency of the proposed methods.
Based on these algorithms, we further explore the self-assembly behavior of liquid crystalline polymers through simulations in four, five, and six dimensions, uncovering intricate three-dimensional spatial structures.
\end{abstract}

%
%
%
%
%
%


\end{frontmatter}



\section{Introduction}
\label{sec:intrd}
The liquid crystalline polymeric systems have attracted much attention due to their potential
industrial applications dependent on customized microstructures\,\cite{pierre1993the}. 
Theoretical studies allow us better understand and predict the phase
behavior of liquid crystalline polymers\,\cite{Schmid1998Schmid, fredrickson2006equilibrium}.
Among various theories, the self-consistent field theory (SCFT)
has been proven to be a powerful tool in the study of phase behaviors of polymeric systems\,\cite{jiang2010Isotropic, Yao2018Topological}. 
Many SCFT studies mainly focus on flexible polymers modelin by continuous
Gaussian chain and achieve great success\,\cite{Masten1994stable, Drolet1999Combinatorial}. 
However, the rigidity of liquid crystal polymer chains heavily affects the phase behaviors of polymeric
systems\,\cite{Matsen1996Melts, Li2014phase, jiang2013self, song2011phase, he2024theory, liu2018archimedean, song2009new}. 
Therefore, a thorough investigation of numerical algorithms for SCFT of liquid crystalline polymeric systems is essential. 

Liquid crystalline polymeric systems typically include semiflexible subchains, which could be modeled by the wormlike chain model. The wormlike chain model\,\cite{saito1967the} can effectively describe the orientation and density distribution of semiflexible subchains, which reults in the emergence of 6-dimensional (6D) partial differential equations (PDEs). Meanwhile, the commonly used long-range interaction in liquid crystals, such as the Maier-Saupe interaction, is introduced into the semiflexible chains to promote parallel alignment between rigid blocks. However, such orientational interactions may lead to a singularity phenomenon, causing a slowdown in the SCFT iteration process. Jiang et al.\,\cite{jiang2013self,jiang2013influence} established a theoretical basis and numerical scheme for studying the semiflexible diblock copolymers. 
They employed a third-order backward differentiation formula in the time direction. However, they did not consider the Maier-Saupe interaction, resulting in the absence of the liquid crystal order. Considering the Maier-Saupe interaction, Liang et al.\,\cite{liang2015efficient} used Fourier pseudo-spectral method and spherical harmonic expansion to discretization spatial and directional variables in flexi-semiflexible diblock copolymers respectively, and used second-order operator splitting method to discretization time variables. Meanwhile, they proposed two effective iterative schemes to solve the SCFT equation. However, these proposed methods in their work are still insufficient to address the computational challenges arising from the high dimensionality and singularity present in current calculations.

In this work, we establish the SCFT model with Maier-Saupe interaction, using a linear diblock flexible-semiflexible molecules as an example, and design a series of algorithms to address the computational challenges arising from high dimensionality and singularity in the model.
When solving PDEs, discretizing high-dimensional PDEs result in unaffordable computational burden. High-order, high-precision numerical algorithms are an appropriate method to alleviate this issue.
For spatial functions, under periodic boundary conditions, we employ Fourier pseudo-spectral methods with exponential convergence\,\cite{ fredrickson2006equilibrium,  jiang2013self, jiang2013influence, Jiang2011Dependence}. For orientation functions, we use spherical harmonic transformations, which also exhibit exponential convergence\,\cite{fredrickson2006equilibrium, jiang2010Isotropic, jiang2013self, song2011phase, jiang2013influence}. Various time discretization schemes can be applied, such as Runge-Kutta type schemes, \,\cite{jiang2013self, jiang2013influence}, operator-splitting scheme\,\cite{fredrickson2006equilibrium, Jiang2011Dependence}. Meanwhile, spectral deferred correction methods \,\cite{ceniceros2019efficient} can be employed to enhance the accuracy of discretization. 
Considering the convection-diffusion characteristics of high-dimensional PDEs, combining finite difference methods with upwind schemes is also an option\,\cite{Yao2018Topological}. 
The effectiveness, accuracy, and computational complexity of these algorithms we proposed are also presented.
To alleviate the slow convergence of SCFT iterations caused by Maier-Saupe interactions, a stable and rapidly converging iterative method is required. The design of nonlinear iterative methods for updating field depends on the mathematical properties of SCFT. Gradient-based iterative methods can be proposed to update field functions when the descent and ascent directions of the free energy can be identified, such as explicit Euler schemes\,\cite{ceniceros2004numerical, jiang2015analytic}, semi-implicit
scheme\,\cite{ceniceros2004numerical, jiang2015analytic, chantawansri2007self}, 
and hybrid conjugate gradient method\,\cite{jiang2015analytic, liang2015efficient}. Anderson mixing methods can be employed to accelerate iterations if SCFT is regarded as a fixed-point problem. However, there is currently no definitive standard for evaluating the efficiency of these iterative methods in SCFT iterations.
In this work, we develop adaptive Anderson mixing methods and cascadic multi-level methods to accelerate SCFT iterations. Moreover, as each ordered structure possesses an optimal period, identifying this period facilitates more precise free energy calculations. We propose an efficient method for the automatic optimization of the computational domain to achieve this.

Our contribution in this paper include: a) developing a suite of highly accurate and efficient algorithms discerning the energy differences among self-assembled structures, b) establishing the evaluation criteria for selecting appropriate high-dimensional PDE solvers based on approximation accuracy, computational complexity and numerical performance the exploration, analysis and implementation of ten numerical methods, c) introducing an adaptive Anderson mixing method and incorporating the cascadic multi-level 
method to accelerate the SCFT calculation, d) designing an efficient method for the automatic optimization of the computational domain, thereby enabling more precise free energy calculations of ordered structures, and e) finding complex three-dimensional spatial structures of liquid-crystalline polymers.

The present paper is organized as follows. In Sec.\,\ref{sec:scft}, we present
a brief description of the SCFT model for the flexible-semiflexible diblock
copolymers. 
In Sec.\,\ref{sec:nume}, we introduce ten time discretization schemes to solve high-dimensional semiflexible propagators. Meanwhile, we present the adaptive Anderson mixing iteration approach and cascadic multi-level method to accelerate the SCFT iteration. Furthermore, we give an approach to automatically optimize the computational box of ordered structures during the iteration process. In Sec.\,\ref{sec:results},  we present numerical results to demonstrate the effectiveness of all proposed algorithm. In Sec.\,\ref{sec:concl}, we provide a comprehensive summary of the study.

\section{Self-consistent field theory}
\label{sec:scft}
\begin{figure}[h]
	\centering
	{
		\includegraphics[width=8cm]{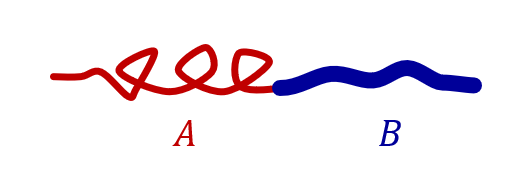}
	}
	\caption{Schematic of the  flexible-semiflexible chain containing a flexible block $A$ (red) and a semiflexible block $B$ (blue).}
	\label{fig:model}
\end{figure}

We conside an incompressible melt consisting of $n$ flexible-semiflexible diblock copolymers in a volume of $V$. 
Each flexible-semiflexible diblock copolymer, with a degree of polymerization $N$, consists of two chemically distinct monomers $A$ and $B$, as schematically shown in Fig.\,\ref{fig:model}.
The number of monomers of the $\alpha$ block  ($i\in  \{A,\, B\}$) is denoted by $N_{\alpha}=f_{\alpha}N$, where $f_{\alpha}$ is the volume fraction. It is noted that $f_{A}+f_{B}=1$, $N_{A}+N_{B}=N$. 
The statistical segment lengths of monomers are $b_{\alpha}$. 
We employ the Gaussian chain model and the wormlike chain model to describe
flexible and semiflexible blocks, respectively\,\cite{fredrickson2006equilibrium}.
The interaction between the monomers $A$ and $B$ is measured by a Flory-Huggins parameter $\chi$, whereas the orientation interaction between semiflexible blocks is depicted by the Maier-Saupe parameter
$\eta$. 
Within the SCFT framework\,\cite{fredrickson2006equilibrium, song2011phase}, the free energy $H$ of flexible-semiflexible diblock copolymer is
\begin{align}
	\begin{split}
		\frac{H}{nk_{B}T}= \frac{1}{V}\int_{V} \left\{\frac{1}{\chi N} \mu_{-}^2(\br) -\mu_+(\br) + \frac{1}{2\eta N} \bM(\br):\bM(\br) \right\}\, d\br -\log Q,
	\end{split}
	\label{H}
\end{align}
where $k_{B}$ is the Boltzmann constant, $T$ is the temperature,
$Q$ is the single chain partition function,
$\bM\in\bbR^{3\times 3}$ is the tensor orientational field. 
The pressure-like field $\mu_{+}$ corresponds to the local incompressibility.
The exchange chemical field $\mu_{-}$ is conjugate to the local density
difference between the polymer densities. 
The mean fields $\omega_{A} (\br)$ and $\omega_{B} (\br)$ are  related to  $\mu_{+}(\br),\,\mu_{-}(\br)$,
\begin{align}
	\omega_A(\br)=\mu_{+}(\br)-\mu_{-}(\br), \quad \omega_B(\br)=\mu_{+}(\br)+\mu_{-}(\br).
	\label{SCFT:w}
\end{align}
The single chain partition function $Q$ can be calculated as
\begin{align}
	\begin{split} 
		Q&=\frac{1}{V}\int_{V}  \int_{\bbS} q_{A}(\br,f) q^{\dagger}_{B}(\br,\bu,1-f) ~d\bu~d\br,\\
		&=\frac{1}{V}\int_{V} \int_{\bbS} q_{B}(\br,\bu,s)q^{\dagger}_B(\br,\bu,1-s)~d\bu~d\br, \quad \forall s\in (f,1],\\
		&= \frac{1}{V} \int_{V} q_{A}(\br,s)q^{\dagger}_{A}(\br,1-s)~d\br,\quad  \forall s\in [0,f], \\
		\label{single:partitionQ}
	\end{split}
\end{align}
where $\bu$ represents the local orientation of the semiflexible segment defined on a unit spherical surface $ \bbS$.
Flexible forward propagator $q_{A}(\br,s)$ describes the 
probability of finding the $s$-th ($0\le s\le f$) segment at spatial position 
$\br$ under the mean field $\omega_{A}(\br)$. 
$q_{A}(\br,s)$ satisfies modified diffusion equation
\begin{align}
	\begin{split}
		\frac{\partial}{\partial s}q_{A}(\br,s)=\,&\nabla^{2}_{\br}q_{A}(\br,s)-\omega_A(\br)q_{A}(\br,s),\quad s\in [0,f],
		\label{propagator:qA:PDE}
	\end{split}
\end{align}
with initial condition $q_{A}(\br, 0)= 1$.

Semiflexible forward propagator $q_{B}(\br,\bu, s)$ describes the 
probability of finding the $s$-th segment at spatial position 
$\br$ with orientation $\bu$ ranging from $s=1$ to $s=1-f$ under the mean field $\omega_{B}(\br)$.
$q_{B}(\br,\bu,s)$ satisfies the ``convection diffusion'' equation
\begin{align}
	\begin{split}
		\frac{\partial}{\partial s}q_{B}(\br,\bu,s)= &-\beta\bu\cdot\nabla_{\br}q_{B}(\br,\bu,s)-\biggl(\omega_{B}(\br)-\bM(\br):\biggr[\bu\bu-\frac{1}{3}\bI\biggm]\biggr)q_{B}(\br,\bu,s)\\&
		+\frac{1}{2\lambda}\nabla_{\bu}^2q_{B}(\br,\bu,s), \quad s\in (f,1],
		\label{propagator:qB:PDE}
	\end{split}
\end{align}
with initial condition $q_{B}(\br,\bu,f)= q_{A}(\br,f)/ 4\pi$. $\beta=\sqrt{6}/v$, $v=b_{A}/(b_{B}\sqrt{N})$  measures the asymmetry of
two monomers, $\lambda$ is the stiffness parameter.
Similarly, the backward semiflexible propagator $q^{\dag}_{B}(\br,\bu,s)$ satisfies 
\begin{align}
	\begin{split}
		\frac{\partial}{\partial s}q^{\dag}_{B}(\br,\bu,s) =\, &\beta\bu\cdot\nabla_{\br}q^{\dag}_{B}(\br,\bu,s)-\biggl(\omega_{B}(\br)-\bM(\br):\biggr[\bu\bu-\frac{1}{3}\bI\biggm]\biggr)q^{\dag}_{B}(\br,\bu,s)\\&
		+\frac{1}{2\lambda}\nabla_{\bu}^2q^{\dag}_{B}(\br,\bu,s), \quad s\in [0,1-f],
		\label{inverse:qB:PDE}
	\end{split}
\end{align}
with initial condition $q^{\dagger}_{B}(\br,\bu,0)=  1/4\pi$.
Backward flexible propagator $q^{\dag}_{A}(\br,s)$ satisfies
\begin{align}
	\begin{split}
		\frac{\partial}{\partial s}q^{\dagger}_{A}(\br,s)=\,&\nabla^{2}_{\br}q^{\dagger}_{A}(\br,s)-\omega_{A}(\br)q^{\dagger}_{A}(\br,s), \quad s\in (1-f,1],
		\label{inverse:qA:PDE}
	\end{split}
\end{align}
with initial condition $q^{\dag}_{A}(\br,1-f)=\int_{\bbS} q^{\dagger}_{B}(\br,\bu,1-f)~d\bu.$

The SCFT equations obtained from the first-order variational derivative of the free energy with respect to field functions are

\begin{align}
	\begin{split}
		&\phi_{A}(\br)=\frac{1}{Q} \int_{0}^{f}q_{A}(\br,s)q^{\dagger}_{A}(\br,1-s)~ds, \\
		&\phi_{B}(\br)= \frac{4\pi}{Q}\int_{f}^{1}\int_{\bbS}  q_{B}(\br,\bu,s)q^{\dagger}_{B}(\br,\bu,1-s)~ d\bu\, ds,\\
		&\bS(\br)=\frac{4\pi}{Q}\int_{f}^{1}\int_{\bbS} q_{B}(\br,\bu,s)\biggl(\bu\bu-\frac{1}{3}\bI\biggr)q^{\dagger}_{B}(\br,\bu,1-s)~ d\bu \,ds,\\
		&\phi_{A}(\br)+\phi_B(\br)-1=0,\\
		&\frac{2}{\chi N}\mu_{-}(\br)- \left[\phi_{A}(\br)-\phi_B(\br)\right]=0,\\
		&\frac{1}{\eta N}\bM(\br)-\bS(\br)=0,\\
	\end{split}
\end{align}
where $\phi_{\alpha}(\br)$ ($\alpha\in {A,\,B}$) and $\bS(\br)$ are the density of the $\alpha$-block and the orientational order parameter, respectively.
The standard saddle point iteration of solving SCFT equations includes
\begin{itemize}
	\item Step 1: Give initial fields $\mu_{+}(\br)$, $\mu_{-}(\br)$, $\bM(\br)$ and computational domain $\mathcal{B}$.
	\item Step 2: Calculate propagators $q_{A}(\br,s)$, $q^{\dagger}_{A}(\br,s)$, $q_{B}(\br,\bu,s)$ and $q^{\dagger}_{B}(\br,\bu,s)$ from current fields.
	\item Step 3: Compute the single chain partition function $Q$, order parameters $\phi_{A}(\br)$, $\phi_{B}(\br)$, $\bS(\br) $ and evaluate the free energy $H$.
	\item Step 4: Output the converged results if the error $\xi$ is less than a prescribed value,  otherwise update fields $\mu_{-}(\br)$, $\mu_{+}(\br)$, $\bM(\br)$, go back to step 2,  where $\xi=E_{\mu_+}+E_{\mu_-}+E_{\bM},\,E_{\mu_+}=\left\| \frac{\delta H}{\delta_{\mu_+}} \right\|,\,E_{\mu_-}=\left\| \frac{\delta H}{\delta_{\mu_-}} \right\|$, $E_{\bM}=\sqrt{\big(\sum_{i,j=1}^{3}\left\| \frac{\delta H}{\delta_{\bM_{i,j}}}
		\right\|\big)/9}$, and the norm $\| \bv \|$, $\bv=(v_{1},v_{2},\cdots, v_{n})\in \bbR^{n}$ is denoted by $\| \bv \|=\sqrt{v_1^2+v_2^2+\dots+v_n^2}/n$.
\end{itemize}

\section{Numerical methods}
\label{sec:nume}
In this section, we systematically study the numerical methods of solving
SCFT with liquid crystalline interaction potential based on the
wormlike chain model. These numerical methods mainly include four parts.
Firstly, we introduce ten time discretization schemes for solving high-dimensional partial differential equations \eqref{propagator:qB:PDE} and \eqref{inverse:qB:PDE}, along with corresponding spatial and orientational discretization methods.
We will compare these methods based on computational complexity analysis, numerical accuracy, and numerical performance to recommend the relatively optimal algorithm.
Secondly, we give several iteration approaches to update the field functions,
including the alternative direction iteration (ADI) method, Anderson mixing (AM)
method and the adaptive Anderson mixing (AAM)
method. Thirdly, we introduce a cascadic multi-level (CML) method to accelerate the SCFT
iteration, which has been used in SCFT
calculation  for flexible polymeric
systems by Ceniceros and Fredrickson\,\cite{ceniceros2004numerical}.
Finally, we propose an optimization method to automatically obtain the optimal computational domain.

\subsection{Discretization scheme of high-dimensional PDEs for semiflexible propagators}
\label{sec:numerical}

The primary numerical challenge in SCFT iterations is the solving of semiflexible propagators. 
Flexible propagators satisfy modified diffusion equations with a d-dimensional spatial variable ($d=2,3$) and a 1-dimensional  time variable. This can be efficiently addressed by using the Fourier pseudo-spectral method and a fourth-order backward difference formula scheme\,\cite{jiang2013self,cochran2006stability, jiang2015self}. Solving semiflexible propagators involve discretizing high-dimensional ``convection-diffusion'' equations with a 3-dimensional  spatial, a 2-dimensional  orientational, and a 1-dimensional  time variable. Due to the addition of orientation in semiflexible propagators, compared to flexible propagators, this presents a significant computational challenge. The goal of this section is to explore the most efficient numerical methods for solving semiflexible propagators. 

In this section, we first introduce ten time discretization schemes,
which can be classified into four groups, Runge-Kutta (RK) type schemes, backward difference formula (BDF) type schemes, operator splitting (OS) schemes, and upwind type schemes. Furthermore, we employ the spectral deferred correction (SDC) method\,\cite{dutt2000spectral} for the time variable to enhance the accuracy. Under the periodic boundary conditions, the Fourier pseudo-spectral method with exponential convergence is used to discretize spatial variables. Specifically, the finite difference scheme is employed for spatial discretization when the upwind scheme is used for time discretization. The orientational variable $\bu$ is discretized by the spherical harmonic expansion. 
Before presenting these numerical methods, Eqn.\,\eqref{propagator:qB:PDE} of the semiflexible propagator
$q_{B}(\br,\bu,s)$  can be written as
\begin{align}
	\begin{split}
		\frac{\partial}{\partial s}q_{B}(\br,\bu,s)
		&~=\mathcal{L}_{1}q_{B}(\br,\bu,s)+\mathcal{L}_{2}q_{B}(\br,\bu,s)+\mathcal{L}_{3}q_{B}(\br,\bu,s),\\
		\mathcal{L}_{1}q_{B}(\br,\bu,s) &:= -\beta\bu\cdot\nabla_{\br}q_{B}(\br,\bu,s),\\
		\mathcal{L}_{2}q_{B}(\br,\bu,s)
		&:=-\bigl(w_{B}(\br)-\bM(\br):[\bu\bu-\frac{1}{3}\bI]\bigr)q_{B}(\br,\bu,s),\\
		\mathcal{L}_{3}q_{B}(\br,\bu,s) &:=
		\dfrac{1}{2\lambda}\nabla_{\bu}^2q_{B}(\br,\bu,s).
	\end{split}
	\label{simple:PDE}
\end{align}

To determine the most suitable method among these algorithms, we conduct an analysis of the computational complexity and time discretization accuracy for each method. Here, we present the semi-discrete scheme of the semiflexible propagator, and a fully discrete scheme is exemplified in  \ref{app:full}. Tab.\,\ref{notation} provides symbols and annotations necessary for illustrating numerical methods.
\begin{table}[H]
	\centering
	\caption{Notations of discretization variables}
	\begin{tabular}{c c}
		\hline
		Notation & Explanation\\
		\hline
		$N_{\br}$&the number of discretization nodes of spatial variable $\br$\\
		$N_s$&the number of discretization nodes of time variable $s$\\	
		$N_\theta$&the number of discretization nodes  of angle  $ \theta $  related to orientation $\bu$ \\
		$N_\varphi$&the number of discretization nodes  of angle $ \varphi $  related to orientation $\bu$ \\	
		$\Delta s$& the time $s$ discretization step size\\
		$\Delta h$ &  the space  $\br$  discretization step size\\
		\hline
	\end{tabular}
	\label{notation}
\end{table}
In practice, the Fast Fourier Transformation (FFT) is used in Fourier pseudo-spectral method, while the SPHEREPACK 3.0\, \cite{Adams1999SPHEREPACK} package is used to implement the fast spherical harmonic transformation and its inverse transformation.
The computation complexity of FFT, inverse Fast Fourier
Transformation (IFFT), Discrete Cosine Transform (DCT), inverse Discrete Cosine
Transform (IDCT) and Spherical Harmonic Transform (SHT) are listed in Tab.\,\ref{tab:complexity:symbol}.
The FFT and IFFT algorithms have the same computational complexity,
\textit{i.e.}, $m_{FFT} = m_{IFFT}$. Similarly, $m_{DCT} = m_{IDCT}$, $m_{SHT} =
m_{ISHT}$.

\begin{table}[H]
	\centering
	\caption{The computational complexity of several transformations}
	\begin{tabular}{c c c}
		\hline
		Algorithm &
		Computational complexity &Notation\\
		\hline
		FFT (IFFT) &$O(N_{\br}\log N_{\br})$&$m_{FFT}$ $(m_{IFFT})$\\
		DCT (IDCT)&$O(N_s\log N_s)$&$m_{DCT}$ $(m_{IDCT})$\\
		SHT(ISHT)&$
		O(N_{\theta}^2N_{\varphi}+N_{\theta}N_{\varphi}\log N_{\varphi})
		$&$m_{SHT}$ $(m_{ISHT})$ \\
		\hline
	\end{tabular}
	\label{tab:complexity:symbol}
\end{table}

\subsubsection{Runge-Kutta type methods}

$\bullet$ \textit{Third-order implicit-explicit (IMEX3) method}

The convection diffusion equations\,\eqref{simple:PDE}  can be expressed as 
\begin{align}
	\frac{\partial}{\partial s}q_{B}(\br,\bu,s)=g(q_{B}(\br,\bu,s))+f(q_{B}(\br,\bu,s)),
	\label{ODE}
\end{align}
where $g(q_{B}(\br,\bu,s))= \mathcal{L}_{1}(q_{B}(\br,\bu,s))$ and $f(q_{B}(\br,\bu,s)) =
\mathcal{L}_{2}(q_{B}(\br,\bu,s))+\mathcal{L}_{3}(q_{B}(\br,\bu,s))$. The IMEX3\,\cite{ascher1995implicit} method discretizes $g(q_{B}(\br,\bu,s))$ implicitly, and $f(q_{B}(\br,\bu,s))$ explicitly.

Concretely, for $g$, we use a
third-stage diagonally-implicit Runge-Kutta scheme\,\cite{wanner1996solving} with coefficients
$A=(a_{ij})\in \bbR^{3\times 3}$, $\bb = (b_i)\in \bbR^{3}$, $\bc = (c_i)\in \bbR^{3}$, 
$i,j = 1,2,3$. The specific values of these coefficients are listed below.
\begin{center}
	\begin{minipage}{\textwidth}
		\begin{minipage}[t]{0.35\textwidth}
			\begin{tabular}{c|ccc} 
				$c_{1}$ & $a_{11}$ & 0 &0\\
				$c_{2}$ & $a_{21}$  & $a_{22}$ & 0\\
				$c_{3}$ & $a_{31}$  & $a_{32}$&$a_{33}$ \\
				\hline
				& $ b_1 $ & $ b_2 $ &  $ b_3 $ 
			\end{tabular}
		\end{minipage}
		\begin{minipage}[t]{0.35\textwidth}
			\begin{tabular}{c|ccc}        
				0.4358665  & 0.4358665 & 0 & 0 \\
				0.7179333 & 0.2820667 & 0.4358665 & 0 \\
				1 & 1.2084966 & -0.6443632 & 0.4358665 \\
				\hline
				& 1.2084966 & -0.6443632 & 0.4358665
			\end{tabular}
		\end{minipage}
	\end{minipage}
\end{center}
For $f$, we adopt a fourth-stage
explicit scheme with $\bar{\bc}=\binom{0}{\bc}$ and coefficient matrices
$\bar{A} = (\bar{a}_{ij})\in \bbR^{4\times 4}$, $i, j = 1, 2, 3, 4$, $\bar{\bb}
= (0,\bb)\in \bbR^{4}$, as shown below.
\begin{center}
	\begin{minipage}{\textwidth}
		\begin{minipage}[t]{0.35\textwidth}
			\begin{tabular}{c|cccc}
				0 & $\bar{a}_{11}$ & 0 &0&0\\
				$c_{1}$  & $\bar{a}_{21}$  & $\bar{a}_{22}$ & 0&0\\
				$c_{2}$	 & $\bar{a}_{31}$  & $\bar{a}_{32}$&$\bar{a}_{33}$&0 \\
				$c_{3}$ & $\bar{a}_{41}$  & $\bar{a}_{42}$&$\bar{a}_{43}$&$\bar{a}_{44}$ \\
				\hline
				&0& $ b_1 $ & $ b_2 $ & $ b_3 $
			\end{tabular}
		\end{minipage}
		\begin{minipage}[t]{0.35\textwidth}
			\begin{tabular}{c|cccc}
				0 & 0 & 0 & 0 & 0 \\
				0.4358665 & 0.4358665 & 0 & 0 & 0\\
				0.7179333 & 0.3212789 & 0.3966544 & 0 & 0\\
				1 & -0.1058583 & 0.5529291 & 0.5529291 & 0 \\
				\hline
				& 0 & 1.2084966 & -0.6443632 &0.4358665
			\end{tabular}
		\end{minipage}
	\end{minipage}
\end{center}
Let $s_{n}$ denote the $n$-th discrete point in time, from $s_{n-1}$ to $s_{n}=s_{n-1}+\Delta s$, applying the IMEX3 method, $q(\br,\bu,s_n)$ can be simplified as $q^n(\br,\bu)$.
\begin{align}
	\begin{split}
		q_{B}^{n}(\br,\bu)&=q_{B}^{n-1}(\br,\bu)+\Delta
		s\sum_{j=1}^{3}b_{j}(K_{j}+\bar{K}_{j+1}).
	\end{split}
\end{align}
Assume that $q_{i}(\br,\bu)$, $i=1,2,3,$ are the uniform interpolation points from $
s_{n-1}$ to $s_n$, then we have
\begin{align}
	K_{i}=g(q_{i}(\br,\bu)),
\end{align}
where
\begin{align}
	q_{i}(\br,\bu)=q_{B}^{n-1}(\br,\bu)+\Delta s\sum_{j=1}^{i}a_{i,j}K_{j}+\Delta
	s\sum_{j=1}^{i}\bar{a}_{i+1,j}\bar{K}_{j}.
\end{align}
and
\begin{align}
	\bar{K}_{1}=f(q_{B}^{n-1}(\br,\bu)),\quad \bar{K}_{i+1}=f(q_{i}(\br,\bu)).
\end{align}
The computational complexity of the IMEX3 method is
\begin{align}
	\begin{split}
		C_{\mbox{IMEX3}}&=8N_sN_{\br}m_{SHT}+9N_{\varphi}N_{\theta}N_sm_{FFT}
		+39N_{\varphi}N_{\br}N_{\theta}N_s\\&+16N_{\varphi}N_{\theta}N_s+
		2N_{\theta}^2 + 2N_{\br}N_{\theta} + 20N_{\varphi}N_{\theta}+ 18N_{\varphi}N_{\br}N_{\theta}.
	\end{split}
\end{align}

$\bullet$ \textit{Fourth-order Runge-Kutta (RK4) method}

The RK4 method\,\cite{wanner1996solving} can also be used to solve the 
Eqn.\eqref{simple:PDE}. Firstly, we rewrite the equation\,\eqref{simple:PDE} into the following form,
\begin{align}
	\frac{\partial}{\partial s}q_{B}(\br,\bu,s)=F(q_{B}(\br,\bu,s)),
	\label{RK:ODE}
\end{align}
where $F(q_{B}(\br,\bu,s))=\mathcal{L}_{1}q_{B}(\br,\bu,s)+\mathcal{L}_{2}q_{B}(\br,\bu,s)+\mathcal{L}_{3}q_{B}(\br,\bu,s)$.
Then, from $s_{n-1}$ to $s_{n}=s_{n-1}+\Delta s$, we denote
\begin{equation}
	\begin{aligned}
		K_{1}&=F(q_{B}^{n-1}(\br,\bu)),
		\\
		K_{2}&=F(q_{B}^{n-1}(\br,\bu)+\Delta sK_{1}/2),
		\\
		K_{3}&=F(q_{B}^{n-1}(\br,\bu)+\Delta sK_{2}/2),
		\\
		K_{4}&=F(q_{B}^{n-1}(\br,\bu)+\Delta sK_{3}).
	\end{aligned}
\end{equation}
Finally, we obtain
\begin{align}
	q_{B}^{n}(\br,\bu)=q_{B}^{n-1}(\br,\bu)+\frac{\Delta s}{6}(K_{1}+2K_{2}+2K_{3}+K_{4}).
\end{align}
The computational complexity of the RK4 method is
\begin{equation}
	\begin{aligned}
		C_{\mbox{RK4}} &=
		8N_sN_{\br}m_{SHT}+8N_{\varphi}N_{\theta}N_sm_{FFT}+25N_{\varphi}N_{\br}N_{\theta}N_s\\
		&+ 2N_{\theta}^2 + 2N_{\br}N_{\theta} + 20N_{\varphi}N_{\theta}+
		26N_{\varphi}N_{\br}N_{\theta}.
	\end{aligned}
\end{equation}
\subsubsection{Backward difference formula (BDF) type methods}

The BDF-type approaches, as one of the linear multistep methods, have been commonly utilized in SCFT calculations\,\cite{jiang2013self, cochran2006stability, jiang2015self}. Here, we also apply this method to solve Eqn.\,\eqref{simple:PDE}.

$\bullet$ \textit{Third-order BDF (BDF3) scheme} 

The BDF3 method has been used to solve semiflexible propagators\,\cite{jiang2013self}, which can be expressed as
\begin{align}
	\begin{split}
		&\frac{11}{6}q_{B}^{n}(\br,\bu)-3q_{B}^{n-1}(\br,\bu)+\frac{3}{2}q_{B}^{n-2}(\br,\bu)-\frac{1}{3}q_{B}^{n-3}(\br,\bu)\\=&\Delta s \mathcal{L}_{1}q_{B}^{n}(\br,\bu)+\Delta s\biggm(3f^{n-1}(\br,\bu)-3f^{n-2}(\br,\bu)+f^{n-3}(\br,\bu)\biggr),
	\end{split}
\end{align}
where $f(\br,\bu,s)=\mathcal{L}_{2}q_{B}(\br,\bu,s)+\mathcal{L}_{3}q_{B}(\br,\bu,s)$.
In our implementation, initial values of the first three steps are obtained through the 
IMEX3 scheme.
The computational complexity of the BDF3 method is
\begin{equation}
	\begin{aligned}
		C_{\mbox{BDF3}} &=
		(2N_s+16)N_{\br}m_{SHT}+(2N_s+12)N_{\varphi}N_{\theta}m_{FFT}+14N_{\varphi}N_{\br}N_{\theta}N_s\\
		&+ N_{\theta}N_{\br}N_s - N_{\br}N_{\theta} +
		20N_{\varphi}N_{\theta}+ 60N_{\varphi}N_{\br}N_{\theta}.
	\end{aligned}
\end{equation}

$\bullet$ \textit{Fourth-order BDF (BDF4) scheme}

The BDF4 method used to solve the Gaussian chain model\,\cite{cochran2006stability} can also be applied to solve wormlike chain model, as illustrated below.
\begin{align}
	\begin{split}
		\frac{25}{12}&q_{B}^{n}(\br,\bu)-4q_{B}^{n-1}(\br,\bu)+3q_{B}^{n-2}(\br,\bu)-\frac{4}{3}q_{B}^{n-3}(\br,\bu)+\frac{1}{4}q_{B}^{n-4}(\br,\bu)\\=&\Delta s\mathcal{L}_{1}q_{B}^{n}(\br,\bu)+\Delta s\biggm(4f^{n-1}(\br,\bu)-6f^{n-2}(\br,\bu)+4f^{n-3}(\br,\bu)-f^{n-4}(\br,\bu)\biggr),
	\end{split}
\end{align}
where $f(\br,\bu,s)=\mathcal{L}_{2}q_{B}(\br,\bu,s)+\mathcal{L}_{3}q_{B}(\br,\bu,s)$.
Initial values of the first four steps can be  obtained through the above-mentioned RK4 method. The computational complexity of the BDF4 method is
\begin{equation}
	\begin{aligned}
		C_{\mbox{BDF4}} &=
		(2N_s+24)N_{\br}m_{SHT}+N_{\varphi}N_{\theta}(2N_s+16)m_{FFT}+18N_{\varphi}N_{\br}N_{\theta}N_s\\
		&+ 2N_{\theta}^2+ 8N_{\varphi}N_{\br}N_s + 30N_{\br}N_{\theta}
		+ 20N_{\varphi}N_{\theta}+ 53N_{\varphi}N_{\br}N_{\theta}.
	\end{aligned}
\end{equation}

\subsubsection{Operator splitting (OS) method}
The second-order operator splitting (OS2) method is a widely used technique for solving Gaussian propagator equations\,\cite{Jiang2011Dependence, gao2013non}. For the semiflexible propagator equation\,\eqref{simple:PDE}, from $s_{n-1}$ to $s_{n}=s_{n-1}+\Delta s$, the OS2 method can be expressed as\,\cite{fredrickson2006equilibrium},
\begin{align}
	\begin{split}
		q_{B}^{n}(\br,\bu) = e^{\mathcal{L}_{1}\Delta s/2}e^{\mathcal{L}_{2}\Delta s/2}
		e^{\mathcal{L}_{3}\Delta s}e^{\mathcal{L}_{2}\Delta s/2}e^{\mathcal{L}_{1}\Delta s/2}
		q_B^{n-1}(\br,\bu).
	\end{split}
\end{align}
The specific steps are as follows:
\begin{description}
	\item [Step 1] Take FFT on $q_B^{n-1}(\br,\bu)$ to obtain $\hat{q}_B^{n-1}(\bk,\bu)$, $\hat{q}_B^{n}(\bk,\bu)=e^{\mathcal{L}_{1}\Delta s/2}\hat{q}_B^{n-1}(\bk,\bu)$,
	\item [Step 2]  Perform IFFT on $\hat{q}_B^{n}(\bk,\bu)$ to compute $q_B^{n}(\br,\bu)$, $q_B^{n}(\br,\bu)=e^{\mathcal{L}_{2}\Delta s/2}q_B^{n}(\br,\bu)$,
	\item [Step 3] Apply SHT on
	$q_B^{n}(\br,\bu)$ to calculate $\hat{q}_{B_{lm}}^{n}(\br)$, $\hat{q}_{B_{lm}}^{n}=e^{\mathcal{L}_{3}\Delta s}\hat{q}_{B_{lm}}^{n}(\br)$, 
	\item [Step 4] Use ISHT on $\hat{q}_{B_{lm}}^{n}(\br)$ to get $q_{B}^{n}(\br,\bu)$, $q_B^{n}(\br,\bu)=e^{\mathcal{L}_{2}\Delta s/2}q_B^{n}(\br,\bu)$,
	\item [Step 5] Execute FFT for $q_B^{n}(\br,\bu)$ to obtain
	$\hat{q}_B^{n}(\bk,\bu)$,
	$\hat{q}_B^{n}(\bk,\bu)=e^{\mathcal{L}_{1}\Delta
		s/2}\hat{q}_B^{n}(\bk,\bu)$, perform IFFT on $\hat{q}_B^{n}(\bk,\bu)$ to obtain $q_B^{n+1}(\br,\bu)$.
\end{description}
The computational complexity of the OS2 method is
\begin{equation}
	\begin{aligned}
		C_{\mbox{OS2}} &=
		2N_sN_{\br}m_{SHT}+4N_{\varphi}N_{\theta}N_sm_{FFT}+5N_{\varphi}N_{\br}N_{\theta}N_s\\
		&+ 2N_{\varphi}N_{\br}N_s + 2N_{\br}N_{\theta} +
		20N_{\varphi}N_{\theta}+ 20N_{\varphi}N_{\br}N_{\theta}.
	\end{aligned}
\end{equation}

\subsubsection{Upwind type schemes}
Due to the nature of the semiflexible propagator equation being a convection-diffusion equation\,\cite{Yao2018Topological}, a common approach to solving such equations is to use upwind schemes to handle the``convection'' term, i.e., the spatial differential term.
Taking 1-dimensional space as an example, explicit and implicit upwind schemes are provided. Similar methods are followed for two-dimensional and three-dimensional spaces.

$\bullet$ \textit{Explicit upwind (EXUP) scheme}

If $\beta\bu\ge 0$,
\begin{align}
	\frac{q^{n+1}_{B}(\br_{j},\bu)-q^{n}_{B}(\br_{j},\bu)}{\Delta s}+\beta\bu\frac{q^{n}_{B}(\br_{j},\bu)-q^{n}_{B}(\br_{j-1},\bu)}{\Delta h}
	-\Gamma(\br_{j},\bu)q^{n}_{B}(\br_{j},\bu)-\frac{1}{2\lambda}\nabla^{2}_{\bu}q^{n}_{B}(\br_{j},\bu)=0,
	\label{EXUP:above}
\end{align}
otherwise $\beta\bu < 0$,
\begin{align}
	\frac{q^{n+1}_{B}(\br_{j},\bu)-q^{n}_{B}(\br_{j},\bu)}{\Delta s}+\beta\bu\frac{q^{n}_{B}(\br_{j+1},\bu)-q^{n}_{B}(\br_{j},\bu)}{\Delta h}
	-\Gamma(\br_{j},\bu)q^{n}_{B}(\br_{j},\bu)-\frac{1}{2\lambda}\nabla^{2}_{\bu}q^{n}_{B}(\br_{j},\bu)=0,
	\label{EXUP:down}
\end{align}
where
$\Gamma(\br,\bu)=-\big(w_{B}(\br)-\bM(\br):\left[\bu\bu-\frac{1}{3}\bI\right]\bigl)$.
The computational complexity of the EXUP method is
\begin{equation}
	\begin{aligned}
		C_{\mbox{EXUP}} &= 
		2N_sN_{\br}m_{SHT}+
		7N_{\varphi}N_{\br}N_{\theta}N_s+2N_{\theta}^2
		+30N_{\varphi}N_{\theta}+20N_{\varphi}N_{\br}N_{\theta}+2N_{\theta}.
	\end{aligned}
\end{equation}

$\bullet$ \textit{Implicit upwind (IMUP) scheme}

If $\beta\bu\ge 0$,
\begin{align}
	\frac{q^{n+1}_{B}(\br_{j},\bu)-q^{n}_{B}(\br_{j},\bu)}{\Delta s}+\beta\bu\frac{q^{n+1}_{B}(\br_{j},\bu)-q^{n+1}_{B}(\br_{j-1},\bu)}{\Delta h}-\Gamma(\br_{j},\bu)q^{n}_{B}(\br_{j},\bu)-\frac{1}{2\lambda}\nabla^{2}_{\bu}q^{n}_{B}(\br_{j},\bu)=0,
	\label{IMUP}
\end{align}
otherwise $\beta\bu < 0$,
\begin{align}
	\frac{q^{n+1}_{B}(\br_{j},\bu)-q^{n}_{B}(\br_{j},\bu)}{\Delta s}+\beta\bu\frac{q^{n+1}_{B}(\br_{j+1},\bu)-q^{n+1}_{B}(\br_{j},\bu)}{\Delta h}-\Gamma(\br_{j},\bu)q^{n}_{B}(\br_{j},\bu)-\frac{1}{2\lambda}\nabla^{2}_{\bu}q^{n}_{B}(\br_{j},\bu)=0.
\end{align}
The computational complexity of the IMUP method is
\begin{equation}
	\begin{aligned}
		C_{\mbox{IMUP}} = 
		2N_sN_{\br}m_{SHT}+20N_{\varphi}N_{\theta}+3N_{\theta}+2N_r+2N_{\theta}^2
		&+
		\frac{11}{3}N_{\varphi}N_{\br}N_{\theta}N_s+18N_{\varphi}N_{\br}N_{\theta}+\frac{1}{3}N_{\varphi}N_{\theta}N_sN_r^3.
	\end{aligned}
\end{equation}

\subsubsection{Spectral deferred correction (SDC) method}
\label{subsubsec:SDC}
The SDC method\,\cite{dutt2000spectral} can be used to improve the time discretization accuracy.
The main idea of the SDC method is to use the spectral
quadrature\,\cite{ceniceros2019efficient}, such as a
Gaussian
or a Chebyshev-node
interpolatory quadrature, to
integrate the time derivative, which can achieve a highly-accuracy numerical solution with a largely
reduced number of quadrature points. The SDC method has
been used to solve flexible propagators\,\cite{ceniceros2019efficient}.
Here we extend it to solve wormlike chain model. 
Before solving, an appropriate method is chosen to obtain an initial solution $ q^{[0]}_B(\br,\bu,s)$.

Firstly,  the exact solution of semiflexible propagator can be given by integrating Eqn.\,\eqref{propagator:qB:PDE} along the time variable $s$,
\begin{equation}
	\begin{aligned}
		q_{B}(\br,\bu,s)= ~&q_{B}(\br,\bu,f)+\int_{f}^{s}\biggl[-\beta\bu\cdot\nabla_{\br}q_{B}(\br,\bu,\tau)+\frac{1}{2\lambda}\nabla^2_{\bu}q_{B}(\br,\bu,\tau)\\&-\Bigl(w_{B}(\br)-\bM(\br):\Bigr[\bu\bu-\frac{1}{3}\bI\Bigm]\Bigr)q_{B}(\br,\bu,\tau)\biggr]
		~d\tau.
	\end{aligned}
	\label{eq:SDC}
\end{equation}
The error between the numerical solution $ q^{[0]}_B(\br,\bu,s)$ and the exact
solution $ q_B(\br,\bu,s) $ is defined as
\begin{align}
	\delta^{[0]}(\br,\bu,s)=q_{B}(\br,\bu,s)-q^{[0]}_{B}(\br,\bu,s),
\end{align}
which satisfies the integral equation 
\begin{align}
	\begin{split}
		\delta^{[0]}(\br,\bu,s)=&\int_{f}^{s}\Bigr[-\beta\bu\cdot\nabla_{\br}\delta^{[0]}(\br,\bu,\tau) +\frac{1}{2\lambda}\nabla^2_{\bu}\delta^{[0]}(\br,\bu,\tau)-\Bigl(w(\br)\\&-\bM(\br):\Bigr[\bu\bu-\frac{1}{3}\bI\Bigm]\Bigr)\delta^{[0]}(\br,\bu,\tau)\Bigm]~d\tau+\epsilon^{[0]}(\br,\bu,s).
		\label{SDC:error}
	\end{split}
\end{align}
where the residual 
\begin{align}
	\begin{split}
		\epsilon^{[0]}(\br,\bu,s)&=q_{B}(\br,\bu,f)+\int_{f}^{s}\biggl[-\beta\bu\cdot\nabla_{\br}q^{[0]}_{B}(\br,\bu,\tau)+\frac{1}{2\lambda}\nabla^2_{\bu}q^{[0]}_{B}(\br,\bu,\tau)\\&-\Bigl(w_{B}(\br)-\bM(\br):\Bigr[\bu\bu-\frac{1}{3}\bI\Bigm]\Bigr)q^{[0]}_{B}(\br,\bu,\tau)\biggr]
		~d\tau-q^{[0]}_{B}(\br,\bu,s).
	\end{split}
\end{align}
Then taking the first derivative of Eqn.\,\eqref{SDC:error} with respect to $s$ leads to
\begin{align}
	\begin{split}
		\frac{d \delta^{[0]}(\br,\bu,s)}{ds}=&\Bigr[-\beta\bu\cdot\nabla_{\br}\delta^{[0]}(\br,\bu,s) +\frac{1}{2\lambda}\nabla^2_{\bu}\delta^{[0]}(\br,\bu,s)-\Bigl(w(\br)\\&-\bM(\br):\Bigr[\bu\bu-\frac{1}{3}\bI\Bigm]\Bigr)\delta^{[0]}(\br,\bu,s)\Bigm]+\frac{d \epsilon^{[0]}(\br,\bu,s)}{ds}.
		\label{SDC:error_derivative}
	\end{split}
\end{align}
$\delta^{[0]}$ and $q_B^{[0]}$ are obtained by using the same numerical method.
Then a corrected numerical solution of $q_B$ is 
\begin{align}
	q^{[1]}_{B}(\br,\bu,s)=q^{[0]}_{B}(\br,\bu,s)+\delta^{[0]}(\br,\bu,s).
\end{align} 

Finally, repeating the above processes, one can have $q_B^{[2]},\cdots ,q_B^{[J]}$, $ J $ is the number of deferred corrections. 
The numerical accuracy of the deferred correction solution along the time variable is
\begin{align}
	\left\| q_B-q_B^{[J]}\right\|=O\left((\Delta s)^{m(J+1)}\right),
	\label{eq_accurate}
\end{align}
where $\Delta s$ is the maximum distance of two adjacent Chebyshev-nodes, $m$ is the order of the chosen numerical scheme to solve the Eqns.\,\eqref{propagator:qB:PDE} and \eqref{SDC:error_derivative}.
In the following, we use
the SDC method to improve the upwind type schemes and the OS2 method.

$\bullet$ \textit{Explicit upwind with spectral deferred correction (EXUP+SDC) method} 

We can use the first-order EXUP scheme with 
variable steps $\Delta s_i, i=1 \cdots n$ to solve $ q^{[0]}_B(\br,\bu,s) $ and $\delta^{[0]}(\br,\bu,s)$.  Subsequently, by applying the SDC method for one correction, we obtain the second-order EXUP+SDC method, and the computational complexity is 
\begin{equation}
	\begin{aligned}
		C_{\mbox{EXUP+SDC}}
		&=6N_sN_{\br}m_{SHT}+2N_{\br}N_{\varphi}N_{\theta}m_{DCT}
		+30N_{\varphi}N_{\theta}\\
		& + 33N_{\varphi}N_{\br}N_{\theta}N_s  +20N_{\varphi}N_{\br}N_{\theta}+2N_{\theta}^2+7N_s-1.
	\end{aligned}
\end{equation}

$\bullet$ \textit{Implicit upwind with spectral deferred correction (IMUP+SDC) method} 

Similarly, $q^{[0]}_B(\br,\bu,s)$ and $\delta^{[0]}(\br,\bu,s)$ can be solved using a variable-step first-order IMUP method.  we apply the SDC method for one correction to obtain the second-order solution with IMUP+SDC scheme. The corresponding computational complexity is 
\begin{equation}
	\begin{aligned}
		C_{\mbox{IMUP+SDC}}
		&=6N_sN_{\br}m_{SHT}+2N_{\br}N_{\varphi}N_{\theta}m_{DCT}+\frac{2}{3}N_{\varphi}N_{\theta}N_sN_r^3\\
		&+20N_{\varphi}N_{\theta} + 
		\frac{79}{3}N_{\varphi}N_{\br}N_{\theta}N_s+18N_{\varphi}N_{\br}N_{\theta}+3N_{\theta}+2N_r+2N_{\theta}^2. 
	\end{aligned}
\end{equation}

$\bullet$ \textit{Second-order operator splitting with spectral deferred correction (OS2+SDC) method} 

The OS2 scheme is employed to solve  $q^{[0]}_B(\br, \bu, s)$ and $\delta^{[0]}(\br, \bu, s)$.  Apply the SDC method for one correction, the resulting OS2+SDC method has fourth-order accuracy.
The corresponding computational
complexity is 
\begin{equation}
	\begin{aligned}
		C_{\mbox{OS2+SDC}}
		&=6N_sN_{\br}m_{SHT}+2N_{\br}N_{\varphi}N_{\theta}m_{DCT}+10N_{\varphi}N_{\theta}N_sm_{FFT}\\
		&+15N_{\varphi}N_{\br}N_{\theta}N_s+ 4N_{\varphi}N_{\br}N_s + 2N_{\br}N_{\theta}
		+ 20N_{\varphi}N_{\theta}+ 20N_{\varphi}N_{\br}N_{\theta}+2N_{\theta}.
	\end{aligned}
\end{equation}
Although SDC method can improve the efficiency of existing algorithms, its computational complexity and computation time are significantly increased. Therefore, we have not considered correcting other algorithms with higher computational complexity in this paper.

Finally, we summarize the dominant computational complexities of all methods, as listed in Tab.\,\ref{tab:primary_complexity}. One can sort the computational
complexity of the ten schemes: $C_{\mbox{IMUP+SDC}}> C_{\mbox{IMUP}}
>C_{\mbox{IMEX3}} > C_{\mbox{RK4}} > C_{\mbox{OS2+SDC}} > C_{\mbox{EXUP+SDC}} >
C_{\mbox{BDF4}} > C_{\mbox{BDF3}} > C_{\mbox{OS2}} > C_{\mbox{EXUP}}$. Therefore, the $\mbox{EXUP}$ method has the least computational amount,
the $\mbox{IMUP+SDC}$
method requires the largest computational cost. 
\begin{table}[H]
	\centering
	\caption{The dominant computational complexity of all proposed methods.}
	\begin{tabular}{c c}
		\hline
		Algorithm &The dominant term of the computational complexity\\
		\hline
		EXUP&$ 2N_sN_{\br}m_{SHT}  $\\
		EXUP+SDC&$ 6N_sN_{\br}m_{SHT}  +2N_{\br}N_{\varphi}N_{\theta}m_{DCT} $\\
		IMUP&$ 2N_sN_{\br}m_{SHT}  +\frac{1}{3}N_{\varphi}N_{\theta}N_sN_r^3 $\\	
		IMUP+SDC&$ 6N_sN_{\br}m_{SHT}
		+2N_{\br}N_{\varphi}N_{\theta}m_{DCT}+\frac{2}{3}N_{\varphi}N_{\theta}N_sN_r^3 $\\
		IMEX3&$8N_sN_{\br}m_{SHT} + 9N_{\varphi}N_{\theta}N_sm_{FFT}$ \\
		RK4&$8N_sN_{\br}m_{SHT}+ 8N_{\varphi}N_{\theta}N_sm_{FFT}$ \\	
		BDF3&$ (2N_s+16)N_{\br}m_{SHT} + N_{\varphi}N_{\theta}(2N_s+12)m_{FFT}$ \\	
		BDF4&$ (2N_s+24)N_{\br}m_{SHT}  + N_{\varphi}N_{\theta}(2N_s+16)m_{FFT}$ \\ 
		OS2&$ 2N_sN_{\br}m_{SHT}  + 4N_sN_{\varphi}N_{\theta}m_{FFT}$ \\ 
		OS2+SDC&$
		6N_sN_{\br}m_{SHT}+2N_{\br}N_{\varphi}N_{\theta}m_{DCT}+10N_{\varphi}N_{\theta}N_sm_{FFT}$ \\ 
		\hline
	\end{tabular}
	\label{tab:primary_complexity}
\end{table}

\subsection{Nonlinear iteration methods}
\label{Sec:Nonlinear}
The iterative methods for updating the fields depend on the mathematical structure of SCFT, and the inclusion of liquid crystalline interaction potentials leads to slow convergence of the iterations. To find an efficient and robust algorithm,
we discuss several common nonlinear iteration methods, including the alternative direction iteration method, the Anderson mixing method and adaptive Anderson mixing algorithm.

\subsubsection{Alternative direction iteration (ADI) method}
In the flexible-semiflexible diblock copolymer system, the ordered structure corresponding to a saddle point
the free energy functional can be obtained by maximizing with respect to $\mu_{+}$, and simultaneously minimizing with respect to $\mu_{-}$ and $\bM$. Therefore, we can update $\mu_{+}$ along with the direction of $\frac{\delta H}{\delta\mu_{+}}$, while update $\mu_{-}$ and $\bM$ along
with the directions of $-\frac{\delta H}{\delta \mu_{-}}$ and $-\frac{\delta H}{\delta \bM}$. Concretely, the ADI method\,\cite{Olsen2008} to update fields can be written as

\begin{subequations}
	\begin{align}
		\frac{\mu^{j+1}_{+}-\mu^{j}_{+}}{\Delta t_{1}}&=(\phi^{j}_{A}+\phi^{j}_{B}-1),\\
		\frac{\mu^{j+1}_{-}-\mu^{j}_{-}}{\Delta t_{2}}&=-\bigl(\frac{2}{\chi N}\mu^{j}_{-}-(\phi^{j}_{A}-\phi^{j}_{B})\bigr),\\
		\frac{\bM^{j+1}-\bM^{j}}{\Delta t_{3}}&=-\bigl(\frac{1}{\eta N}\bM^{j}-\bS^{j}\bigr),
	\end{align} 
	\label{method:ADI}
\end{subequations}
where $\Delta t_{i}$, $i=1, 2, 3$ is the iteration step length.

\subsubsection{Anderson mixing (AM) approach}
\label{subsec:anderson}
AM method was originally proposed
by Anderson to solve nonlinear integral equation\,\cite{Anderson1965Iterative}, which makes full use of the historical information of the previous $m$ step to accelerate convergence. Recently, the Anderson mixing method has been widely
used in the polymeric SCFT calculation\,\cite{Arora2017, Pratapa2016, Pollock2019}. It has greatly accelerated convergence in SCFT calculation of
flexible systems. In this work, we try to discuss its performance for the semiflexible polymeric systems with the
liquid-crystal interaction potential. Considering a fixed point $\psi(\bx)=\bx$, the AM method updates the field functions $\bx = \{\mu_{+}, \mu_{-}, \bM\}$ by applying historical information from prior $m$ steps. Algorithm \ref{alg:anderson} presents the implementation of the Anderson mixing method.

\begin{algorithm}[t]
	\caption{ Anderson mixing algorithm}
	\label{alg:anderson}
	\KwIn{$\bx_{0}$, $\psi$, $m$.}
	
	$\bx_{1}=\psi(\bx_{0})$, $F_{0}=\psi(\bx_{0})-\bx_{0}$,
	
	\For{$k=1$, $2$, $\cdots$, $ n-1 $}{
		$m_{k}=\min(m,k)$,
		
		$F_{k}=\psi(\bx_{k})-\bx_{k}$,
		
		Minimize~ $\big\|\sum_{j=0}^{m_{k}}\alpha^{(k)}_{j} F_{k-m_{k}+j}
		\big\|$ subject to $\sum_{j=0}^{m_{k}}\alpha^{(k)}_{j}=1$,
		
		$\bx_{k+1}=(1-\beta_{k})\sum_{j=0}^{m_{k}}\alpha^{(k)}_{j}\bx_{k-m_{k}+j}+\beta_{k}\sum_{j=0}^{m_{k}}\alpha^{(k)}_{j}\psi(\bx_{k-m_{k}+j})$,
		where $\beta_{k}$ is the weight.
	}	
	\KwOut{converged value $\bx_n$.}
\end{algorithm}

To avoid the instability and improve the convergence speed of AM method, we propose an adaptive AM method. The main idea of the  adaptive AM method is to employ the ADI method within the AM framework when an increase in error $\xi$ is detected. Compared with the traditional AM method, this strategy efficiently avoids ineffective iterations and realigns the convergence trajectory towards the optimal path. Details of the adaptive AM are given in Algorithm\,\ref{alg:aanderson}.
\begin{algorithm}[htbp!]
	\caption{ Adaptive Anderson mixing algorithm}
	\label{alg:aanderson}
	\KwIn{$\bx_{0}$, $\psi$, $m$.}
	\For{$k=1$, $2$, $\cdots$, $ n-1 $}{
		carry out Alg.\,\ref{alg:anderson}
		
		\If{$\xi_k > \xi_{k-1}$}{use the ADI method}
		carry out Alg.\,\ref{alg:anderson}
	}
	\KwOut{converged value $\bx_n$.}
\end{algorithm}

\subsection{Cascadic multi-level (CML) method}
\label{mutilevel}

The CML method is an efficient acceleration technique for nonlinear iterations, which was originally introduced to SCFT calculations for flexible polymeric chain systems by Ceniceros and Fredrickson\,\cite{ceniceros2004numerical}.
The idea behind the CML iteration is to rapidly reduce high-frequency errors on fine grids and eliminate low-frequency errors on coarse grids simultaneously through smoothing or relaxation methods. Here we extend the CML strategy to flexible-semiflexible systems.

Keep the symbols in Tab.\,\ref{notation}, and we define $(N_{\br}, N_{\theta}, N_{\varphi}, N_{s})$ as the finest grid level.
One can choose the 
coarsest level $(N_{\br}^0, N_{\theta}^0, N_{\varphi}^0, N_{s}^0)$ such that $(N_{\br}, N_{\theta}, N_{\varphi}, N_{s}) = (2^L N_{\br}^0, 2^L N_{\theta}^0, 2^L N_{\varphi}^0, 2^L N_{s}^0)$, where $L$ is an integer. Then we can obtain $L+1$ level grids:
$$
(N_{\br}^i, N_{\theta}^i, N_{\varphi}^i, N_{s}^i) = (2^iN_{\br}^0, 2^iN_{\theta}^0, 2^iN_{\varphi}^0, 2^iN_{s}^0), \quad i=0,1, \dots L,
$$
with $(N_{\br}^L, N_{\theta}^L, N_{\varphi}^L, N_{s}^L)=(N_{\br}, N_{\theta}, N_{\varphi}, N_{s})$. 
On each level $i$, let $\calW^i$ be the initial fields at a given grid $(N_{\br}^i, N_{\theta}^i, N_{\varphi}^i, N_{s}^i)$, $\xi_i$ denote the error tolerance and $\calW_c^i$ be the convergence fields. 
With these notations, the CML method can be described simply as the algorithm\,\ref{alg:CML}.

\begin{algorithm}[htbp!]
	\caption{Cascadic multi-level method}
	\label{alg:CML}
	
	\KwIn{
		Extract the coarsest grid fields $\calW^{0}$ from the current fields.}
	
	\For{$k=0$, $1$, $2$, $\cdots$, $ L-1 $}{

		\If{$Error > \xi_i$}{continue SCFT iteration.}	
		
		\Else{obtain the  fields $\calW_c^{i}$ on the level $i$.}
		obtain the fields $\calW^{i+1}$ on the
		level $i+1$ by interpolating $\calW_c^{i}$,  		
	}	
	\KwOut{ Convergent fields $\calW_c^{L}$ on the level $L$.}
\end{algorithm}

In algorithm\,\ref{alg:CML}, efficient interpolation involves initially zero-padding the Fourier coefficients of $\calW^{i-1}$ to obtain the Fourier coefficients of $\calW^{i}$ on a finer grid $(2N_{\br}^l, 2N_{\theta}^l, 2N_{\varphi}^l, 2N_{s}^l)$, followed by executing IFFT to obtain $\calW^{i}$.
	
		
		


\subsection{Optimizing computational domain}
\label{calcuad}
The algorithm previously presented focuses on finding the saddle point of the free energy in a fixed computational domain.
Note that Each crystal structure has its own periodicity. Finding the optimal domain for each structure is crucial in SCFT calculations. 
Because the free energy $H$ is functional with respect to the calculation domain, the optimization can be achieved by minimizing the free energy relative to the computational domain. 
Denoting the computational domain as $\mathcal{B}$, the first-order optimal condition of finding the optimum $\mathcal{B}$ is 
\begin{equation}
\left\|\frac{\partial H[\calB]}{\partial \calB} \right\|_{\ell^{\infty}}= 
\max \left|\frac{\partial H[\mathcal{B}]}{\partial b_{ij}}\right|=0, \qquad i,j = 1, \cdots, d.
\end{equation}
In practice, we always choose a proper coordinate system such that $b_{ij} = 0$ when $1 \leq i < j \leq d$. Using this
strategy, the variables are required to be optimized reduce to $d(d + 1)/2$. There are many optimization methods
to solve the unconstrained minimization problem\,\cite{nocedal2006numerical}. In this work, we use the Barzilai-Borwein method\,\cite{Barzilai1988two} to
update the computational box $\mathcal{B}$. When the computational domain is adjusted at $k$-th time, the quantities are marked by the superscript $k$. For each entry $b_{ij}$,

\begin{equation}
b_{ij}^{k+1}=b_{ij}^{k}-\alpha^k\frac{\partial H[\mathcal{B}^{k}]}{\partial b_{ij}^{k}},
\label{adjust_B}
\end{equation}
where $1 \leq i < j \leq d$. Eqns.\,\eqref{adjust_B} can be written
as the matrix form
\begin{equation}
\mathcal{B}^{k+1}=\mathcal{B}^{k}-D^k\frac{\partial H[\mathcal{B}^k]}{\partial
	\mathcal{B}^k},
	\end{equation}
	where $D^k=\alpha^{k}I$. To make $D^k$ have the  quasi-Newton property \cite{Dai2013}, we
	compute $ \alpha^{k} $ such that
	\begin{equation}
\alpha^k = \argmin_{\alpha} \parallel s^{k-1}-D y^{k-1}\parallel_{\ell^\infty},
\label{eq:bb}
\end{equation}
where 
$$
s^{k-1}=\mathcal{B}^{k}-\mathcal{B}^{k-1}  ,  \quad y^{k-1}=\frac{\partial
H[\mathcal{B}^k]}{\partial \mathcal{B}^k}-\frac{\partial
H[\mathcal{B}^{k-1}]}{\partial \mathcal{B}^{k-1}} .
$$
Analytically calculating $\frac{\partial H[\mathcal{B}^{k}]}{\partial \mathcal{B}^{k}}$
is difficult, thus we compute it by the central difference method.
The solution of the Eqn.\,\eqref{eq:bb} can be written as \begin{equation}
\alpha^{k}=\frac{{s^{k-1}}^{T}s^{k-1}}{{s^{k-1}}^{T}y^{k-1}}.
\end{equation}

\begin{figure}[H]
\centering
{
	\includegraphics[width=14cm]{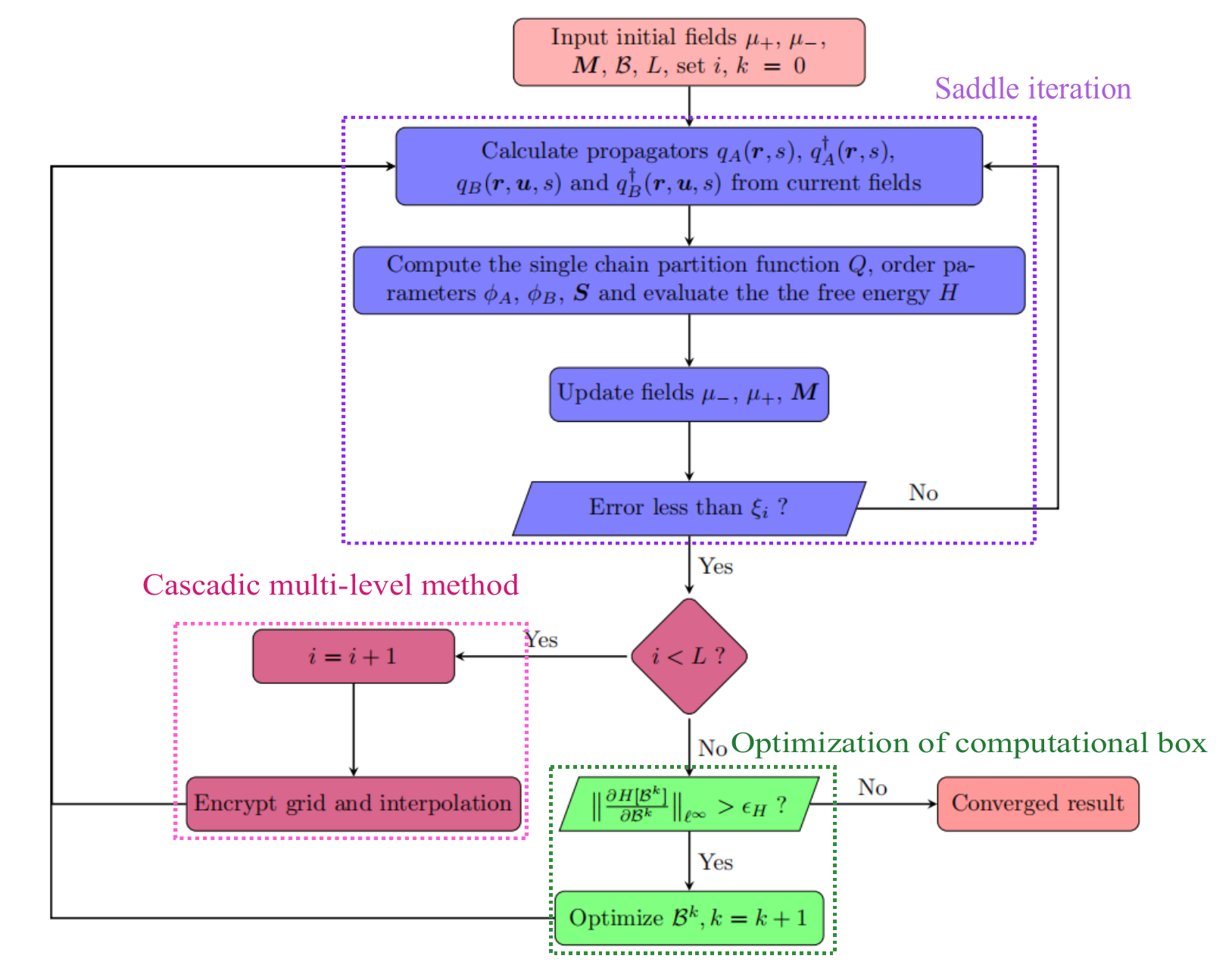}
}
\caption{The flowchart of SCFT iteration process.}
\label{fig:SCFT}
\end{figure}

Note that the proposed method can optimize the size and shape of the computational area automatically during the
optimization process.
The iterative procedure of finding the saddle point of the SCFT equations and optimizing the computational domain is summarized as 
\begin{itemize}
\item [Step 1] Give an initial box $\mathcal{B}^0$ and the initial fields, set $k = 0$.
\item [Step 2] Fix box $\mathcal{B}^k$, perform SCFT iteration until convergence,
calculate the free energy  $H[\mathcal{B}^k]$ and the derivative $\frac{\partial H[\mathcal{B}^{k}]}{\partial \mathcal{B}^{k}}$.
\item [Step 3] If $\left\|\frac{\partial H}{\partial \mathcal{B}^k}\right\|_{\ell^\infty}> \epsilon_H$, compute the step size $\alpha^{k}$ and optimize  $\mathcal{B}^k$ by the Barzilai-Borwein method, set $k = k + 1$, go back
to step 2, otherwise output $H[\mathcal{B}^{k}]$, $ \mathcal{B}^{k} $.
\end{itemize}
The completed SCFT iteration procedure, including saddle point iteration, CML method, and optimal computational domain is presented in Fig.\,\ref{fig:SCFT}.

\section{Numerical results and discussion}
\label{sec:results}

In this section, we present several numerical examples to demonstrate the performance of our proposed methods.
The BDF4 scheme has been used to solve the flexible propagators. All calculations are executed on the CPU through the Intel(R) Xeon(R) Gold 6330
CPU with 28 cores.

\subsection{Efficiency of time discretization schemes for high-dimensional PDEs}
\label{timediscrete}

In this section, all numerical experiments take the case of 4D (space 1D + orientation 2D + time 1D) variables as an example, with parameters listed in Tab.\,\ref{1D:parameter} and initial values referenced from\,\cite{liang2015efficient}.
\begin{table}[H]
\centering
\caption{The parameters on the 4D calculation.}
\begin{tabular}{cccccccccc}
	\hline
	$\chi N$ & $\eta N$ & $f$ & $\calB$ & $ \beta
	$&$N_{\br}$&$N_{\theta}$&$N_{\varphi}$\\
	\hline
	16 & 8  & 0.6 & $ 1.31 $ & 2 & 64 & 8 & 17\\
	\hline
\end{tabular}
\label{1D:parameter}
\end{table}

Fixed spatial and orientational variables are shown in Tab.\,\ref{1D:parameter}, ensuring accuracy in both spatial and orientational dimensions. 
To verify the numerical accuracy of ten time discretization schemes, we set $N_s=20$ and successively refine by bisection, obtaining corresponding errors and ultimately determining the order of accuracy. The error is defined as $\|q^* - q_h\|$, where $q^*$ represents the ``numerically exact solution'' obtained through the RK4 method with 2560 time discretization nodes, and $q_h$ represents the numerical solution. One can find that the numerical errors and convergence orders for the ten time discretization schemes used to solve the semiflexible propagator in Tab.\,\ref{tab:RKandBDF}. All these order are consistent with theoretical results. Obviously, the errors of RK4 and OS2+SDC methods can reach $2 \times 10^{-12}$ with an order of 4, making them the most accurate among all methods.

\begin{table}[H]
\centering
\caption{The error and numerical order of ten numerical discretization schemes for
	solving the semiflexible propagators, with model parameter $\lambda=\infty$.}
\begin{tabular}{ccccccccccc}
	\hline
	\multirow{2}*{$N_s$}&\multicolumn{2}{c}{RK4}&
	\multicolumn{2}{c}{IMEX3}&\multicolumn{2}{c}{BDF4}&\multicolumn{2}{c}{BDF3}&\multicolumn{2}{c}{OS2}\\
	\cmidrule(lr){2-11}
	&  Error& Order & Error& Order & Error& Order & Error & Order&Error& Order \\
	\hline
	20&1.37e-07& - & 1.10e-06 &  - &   7.96e-06&  - &5.64e-05& - &1.21e-04&-\\
	40&8.40e-09&  4.02 &  1.38e-07 & 3.00 &  5.48e-07  & 3.86&   7.43e-06 & 2.92&3.03e-05&1.99\\
	80&5.18e-10&  4.01 & 1.73e-08 & 2.99  & 3.51e-08  & 3.96 &  9.46e-07& 2.97&7.57e-06 &2.00\\
	160&3.21e-11&  4.01 &  2.16e-09 & 2.99 &  2.21e-09&   3.98 &  1.18e-07 & 2.97&1.89e-06&2.00\\
	320&2.00e-12&  4.00 &  2.71e-10 & 2.99 &  1.39e-10&   3.99 &  1.49e-08 & 2.99&4.71e-07&2.01\\
	
	\hline
	\multirow{2}*{$N_s$}&\multicolumn{2}{c}{EXUP}&
	\multicolumn{2}{c}{IMUP}&\multicolumn{2}{c}{EXUP+SDC}&\multicolumn{2}{c}{IMUP+SDC}&\multicolumn{2}{c}{OS2+SDC}\\
	\cmidrule(lr){2-11}
	&  Error& Order & Error& Order & Error& Order & Error & Order& Error & Order\\
	\hline
	20&2.13e-03 & -  & 1.86e-03  & - & 2.13e-04 &  -  & 2.70e-05&-&1.39e-07&-\\
	40&1.03e-03 & 1.04  &  9.19e-04 & 1.02 &   4.99e-05&  2.09   & 7.01e-06 & 1.94&8.90e-09&3.97\\
	80&5.05e-04 & 1.02  &  4.55e-04 & 1.01 &   1.18e-05 & 2.07  &  1.78e-06 & 1.97&5.59e-10&3.99\\
	160&2.49e-04&  1.02 &   2.25e-04&  1.01 &   2.87e-06  &2.04&    4.50e-07 & 1.98&3.50e-11&3.99\\
	320&1.22e-04 & 1.02  &  1.10e-04 & 1.02  &  7.05e-07 & 2.02  &  1.11e-07  &2.00&2.21e-12&3.98\\
	\hline
\end{tabular}
\label{tab:RKandBDF}
\end{table}

\begin{figure}[htpt!]
\centering
{
	\includegraphics[width=6.5cm]{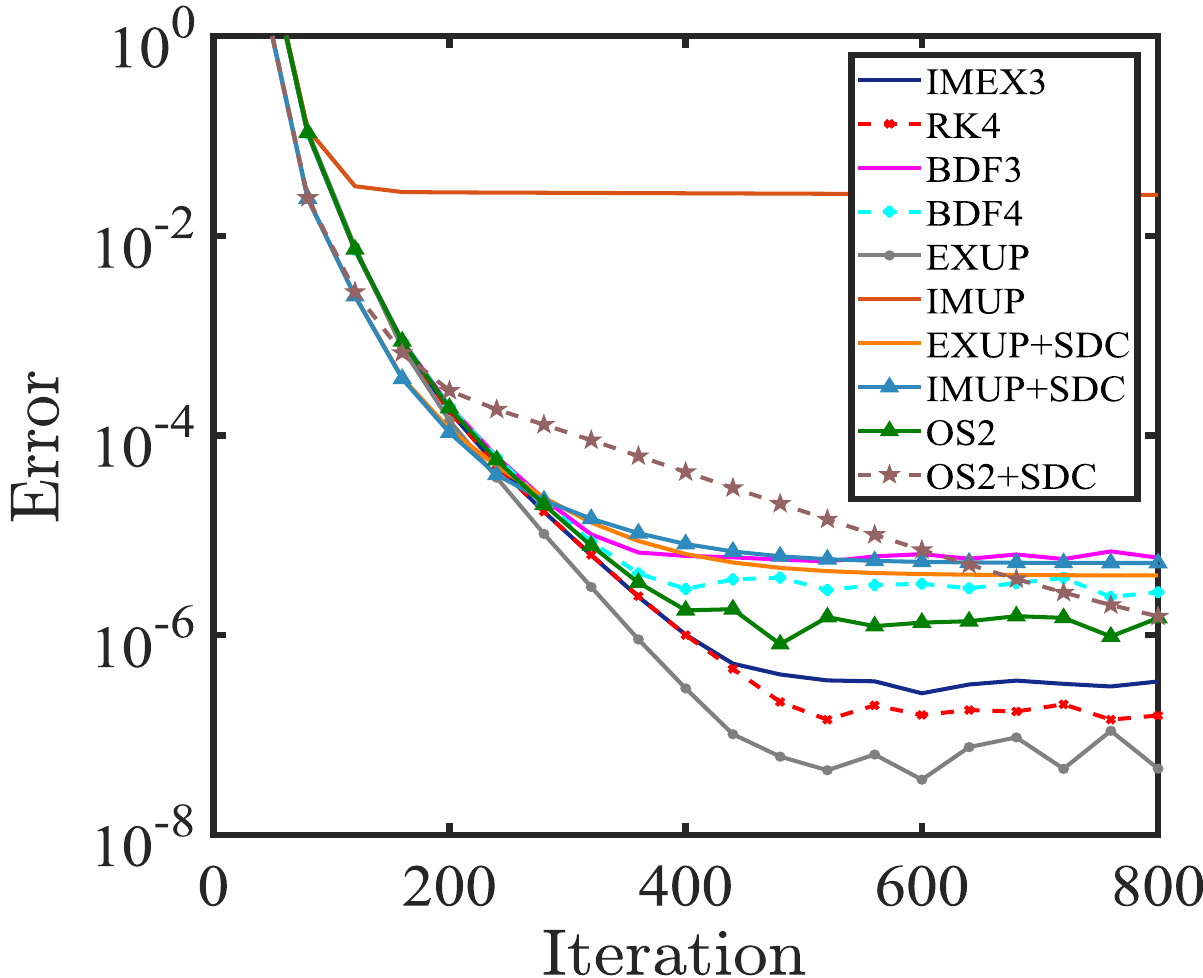}
}
\caption{Comparison of ten time discretization schemes in
	Sec.\,\ref{sec:numerical}, with model
	parameter $\lambda=100$.}
\label{errortotal1D}
\end{figure}

To further compare the performance of ten different time discretization methods, we apply them
to SCFT calculations. The ADI method are chosen as the nonlinear iteration method to guarantee the stability of iteration process. Fig.\,\ref{errortotal1D} shows the iteration errors for solving semiflexible propagators using various time discretization schemes within 800 iteration steps. Notably, the EXUP method reduces the error to approximately $\mathcal{O}\left(10^{-6}\right)$ in just 358 iterations. The RK4 scheme and the IMEX3 method follow closely, both  requiring 402 iterations. Compared to other algorithms, these three methods converge the faster and achieve the lower errors.

\begin{figure}[htpt!]
\centering
{
	\includegraphics[width=6.5cm]{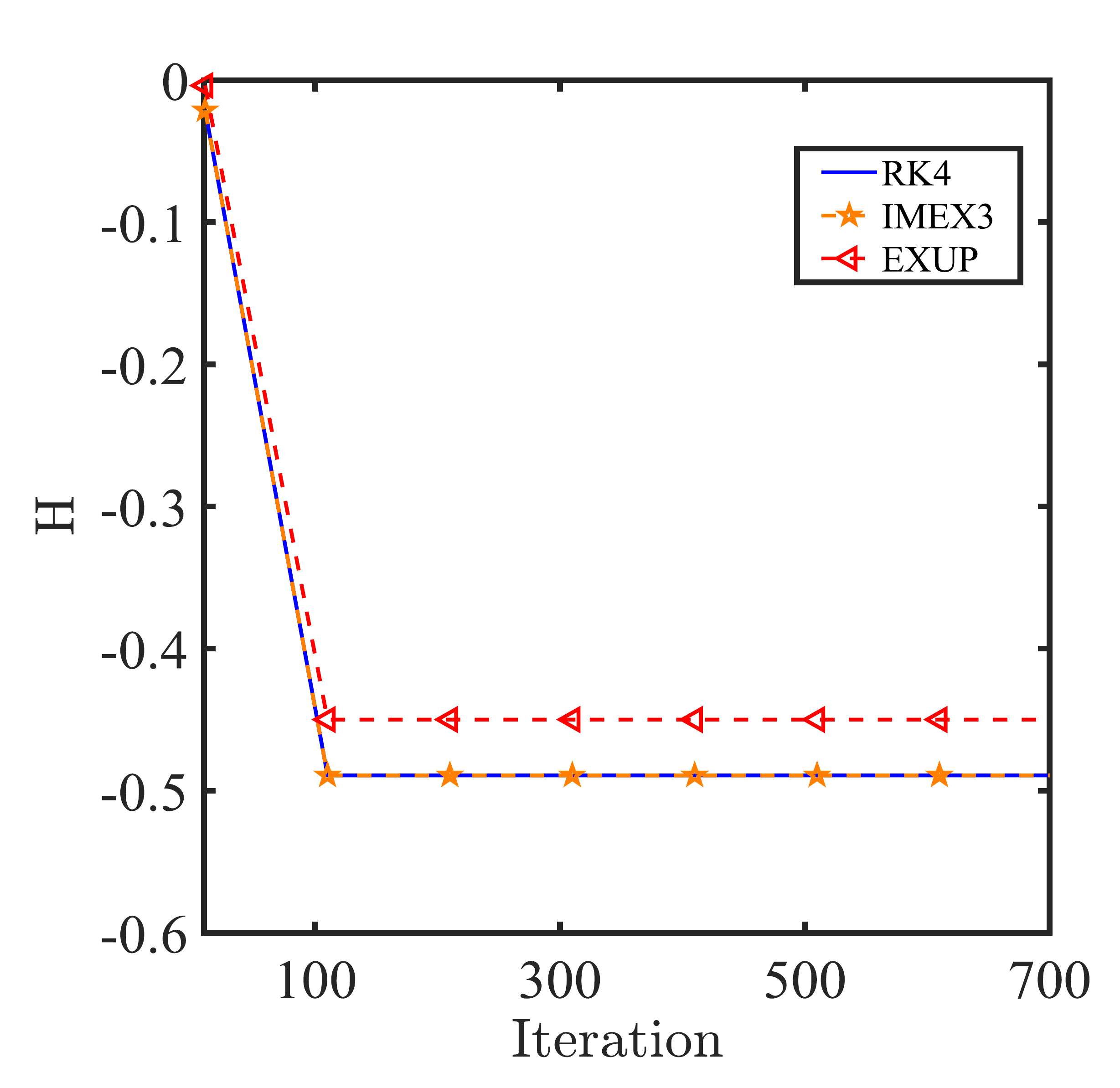}
}
\caption{Numerical behavior of free energy $H$ in the RK4 method, the
	IMEX3 method and the EXUP method, with the parameter $\lambda=100$.}
\label{errortotal1D2}
\end{figure}

To further observe the converged accuracy of the
three methods, Fig.\,\ref{errortotal1D2} expresses that the changes of free energy H. Compared to the EXUP method, the RK4 and IMEX3 methods reach lower free energy $H$ when free energy converges, demonstrating greater stability of these two numerical schemes in the SCFT system. From Tab.\,\ref{tab:primary_complexity} and Tab.\,\ref{tab:RKandBDF}, RK4 method presents higher numerical accuracy and lower computational complexity compared to IMEX3. Therefore, we select RK4 method as the time discretization method for the semiflexible propagator.

\subsection{Iteration efficiency}
\label{subsec:iteration}
In this section, we focus on the effectiveness of nonlinear iterative algorithms. Firstly, we test the performances of these schemes as described in Sec.\,\ref{Sec:Nonlinear} using three iteration methods as follows: (i) ADI scheme, (ii) Anderson mixing method, and (iii) adaptive Anderson mixing method. Secondly,
we use the CML technique to accelerate SCFT
iterations. 
The numerical experiments are based on the 4-dimensional case (space 1D + orientation 2D + contour 1D), with parameters $N_\mathbf{r} = 64$, $N_\theta=16$, $N_{\psi}=33$, $N_s = 200$. No specified, we set $\beta=9.8$, $\eta N = 8$, $\lambda=200$, $f=0.6$, and $\chi N = 16$.

\begin{figure}[H]
\centering
\includegraphics[width=6cm]{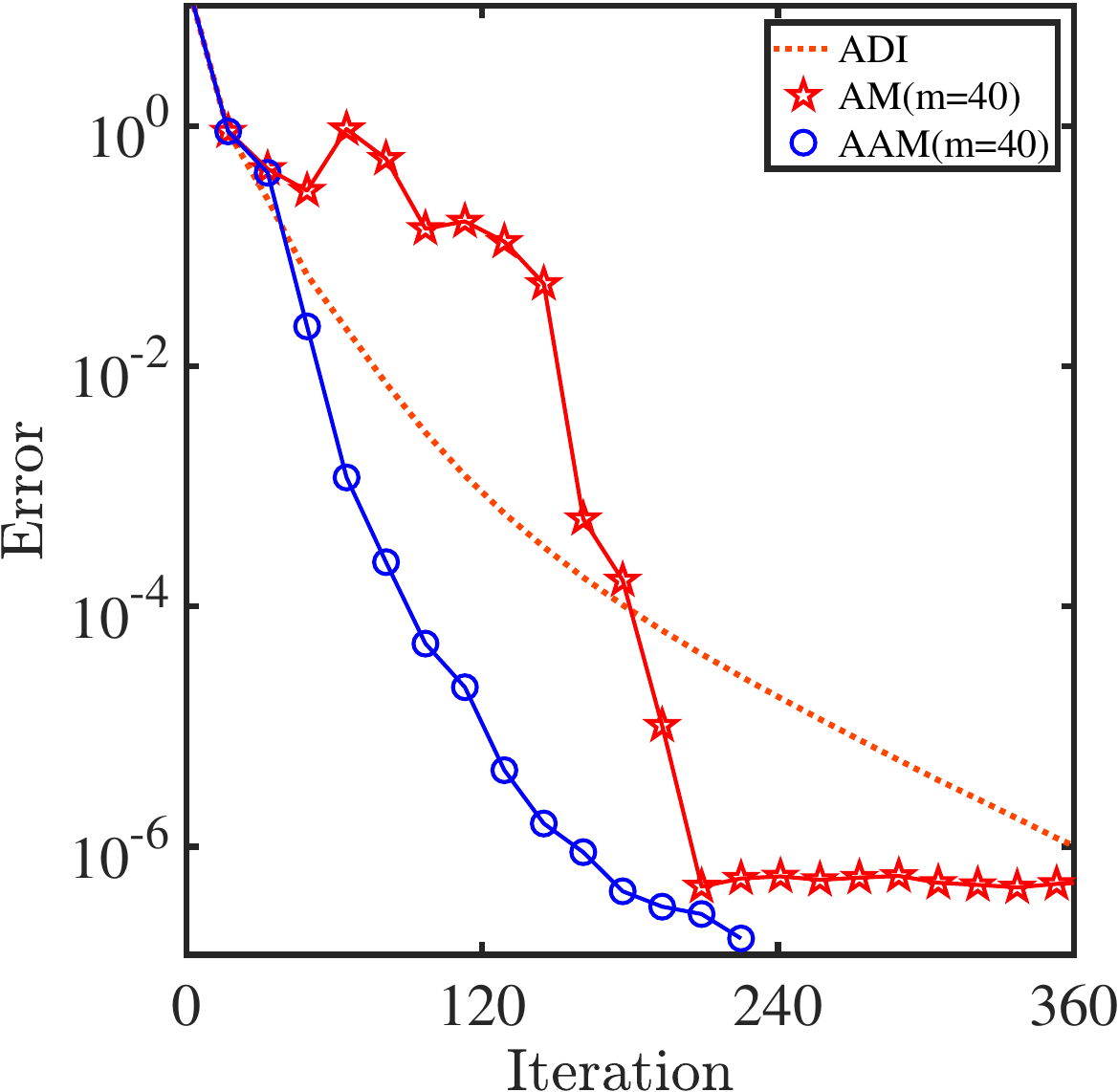}
\caption{Comparison of the ADI method, AM method and AAM
	method, with  parameters $\chi N=16$, $\eta N=8$, $f=0.6$,
	$\lambda=100$, $\beta=2$, $N_{s}=200$, $N_{\br}=32$, $N_{\theta}=8$ and 
	$N_{\varphi}=17$.}
\label{fig:Anderson:compare}
\end{figure}

\begin{figure}[H]                                                               
\centering                                                                  
\includegraphics[width=6.5cm]{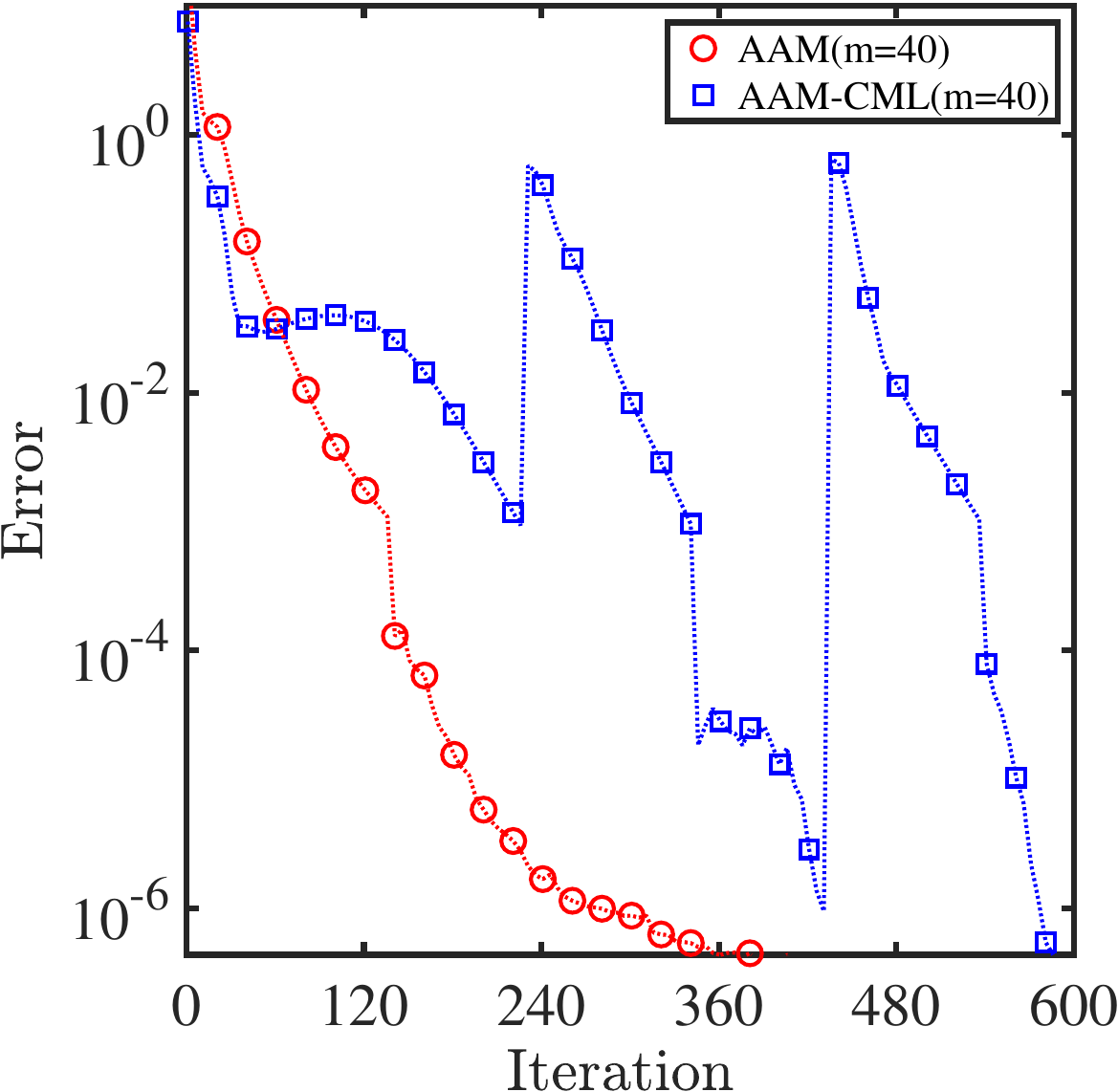}              
\caption{Comparison of the AAM method  and the AAM-CML
	method, with parameters $\chi N=16$, $\eta N=8$, $f=0.6$,
	$\lambda=100$, $\beta=2$, $N_{s}=200$, $N_{\br}=32$, $N_{\theta}=8$ and 
	$N_{\varphi}=17$.}
\label{Anderson:compareaam}
\end{figure}

The convergence of ADI method, AM method and AAM method are shown in Fig.\,\ref{fig:Anderson:compare}.
One can find that three iterative methods can reach the saddle points in the SCFT iterations. Among these methods, the ADI method has the slowest converged rate, the AAM method converges fastest and has fewer iteration steps than the AM method with parameter $m = 40$. Note that the value of $m$ depends on the case. 
Specifically, the AAM method achieves an error of $1.7 \times 10^{-7}$ in only 225 steps, while the AM and ADI methods require 360 steps to achieve errors of $5.3 \times 10^{-7}$ and $1 \times 10^{-6}$, respectively. Therefore, the AAM method is currently the most recommended iterative scheme.
Then we use the CML method to accelerate the AAM method in SCFT iteration. Let the coarsest grid $(N_{\br}^0, N_\phi^0, N_\theta^0, N_s^0)$ be $(16,9,4,50)$ and three levels ($L=3$) are used in this case. The error tolerance $\xi$ in the three grid levels are set as  $1\times 10^{-4}$, $1\times 10^{-6}$ and $3\times10^{-7}$, respectively.  Fig.\,\ref{Anderson:compareaam} shows the SCFT iteration error for the AAM and its CML version (AAM-CML).
At each grid change, there is a significant jump in the error, but it quickly drops to the required error tolerance, eventually achieving the same level of accuracy as AAM method. Although the AAM-CML method takes more steps, its CPU time is significantly reduced, saving $63.88\%$ CPU time  compared to AAM method, as shown in Tab.\,\ref{CPU:XN=6}.

\begin{table}[H]
\caption{Comparison of the CPU time between the AAM and 
	AAM-CML method. }
\centering
\scalebox{0.9}{
	\begin{tabular}{ccc}
		\hline
		AAM ($t$) & AAM-CML ($t_c$) &  $1-t_c/t$~(\%)\\
		\hline
		4710.18s & 1701.12s& $63.88\%$ \\
		\hline
\end{tabular}}
\label{CPU:XN=6}
\end{table}


\subsection{Ordered phases}
\label{subsec:orderedstructures}

In this section, we investigate self-assembled patterns of flexible-semiflexible diblock copolymers to demonstrate the effectiveness of our approaches in 4D (space 1D), 5D (space 2D), 6D (space 3D) (orientation 2D + time 1D) variables, respectively. 
The termination criterion for self-consistent field iterations is the total gradient difference $\xi < 1 \times 10^{-5}$ or the energy difference is less than $1\times 10^{-8}$ between two  sequential iterations.

In the study of flexible-semiflexible diblock copolymers, both nematic and smectic phases have been observed in theoretical\,\cite{masten1998liquid, Semenov1992Theory} and experimental\,\cite{chen1996smectic} researches. To demonstrate the effectiveness of our algorithm, we select the classic nematic and smectic-C phases as examples for
initial values in 4D calculations.

\begin{figure}[H]
\centering
{
	\includegraphics[width=12cm]{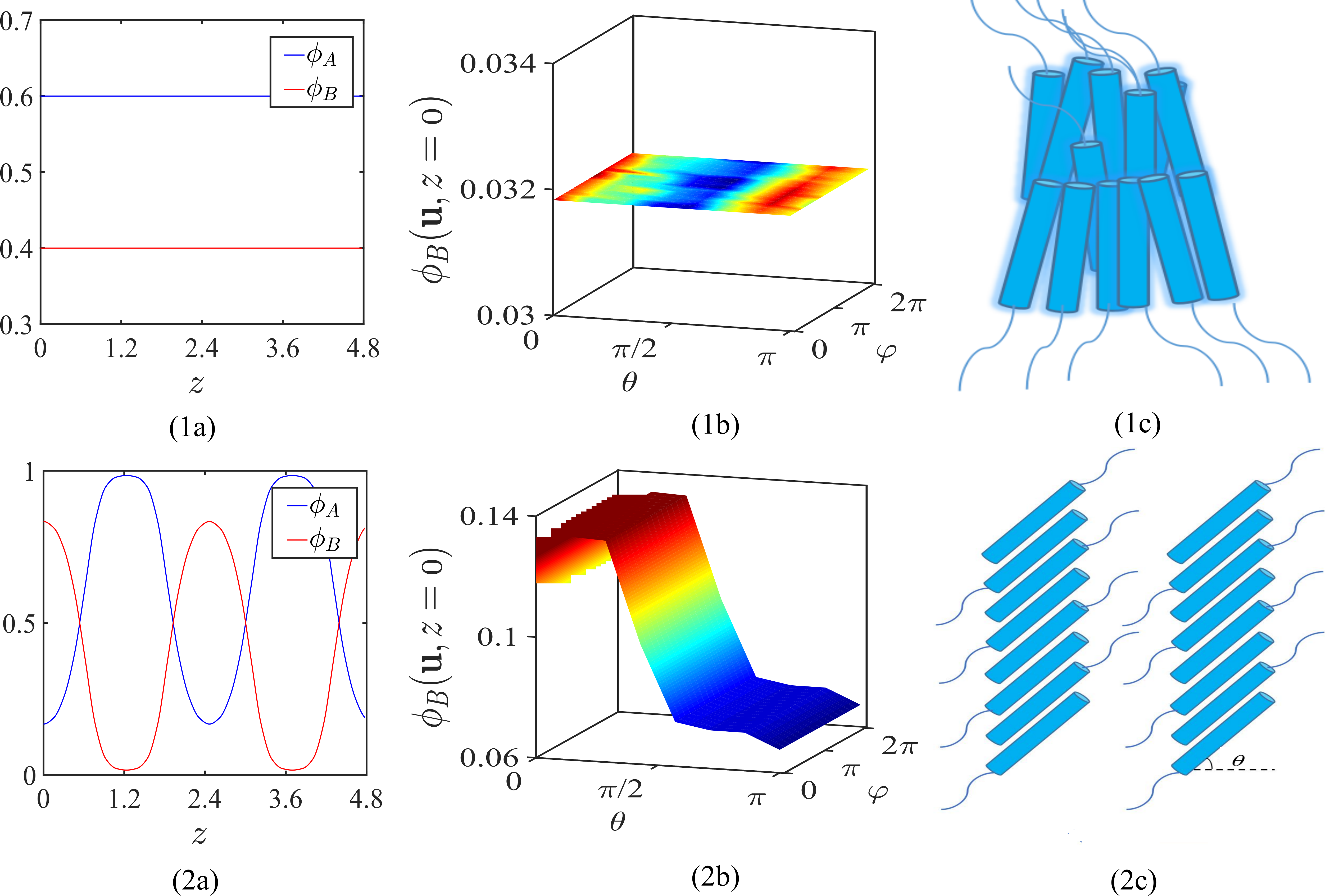}
}
\centering
\caption{(1) Nematic phase, (2) smectic-C phase.
	The first, second, and third columns show the density distributions, orientation distribution, and schematic arrangement diagrams, respectively.}
\label{Fig:nematic}
\end{figure}

\begin{table}[H]
\centering
\caption{The parameters for nematic and smectic-C phases.}
\begin{tabular}{cccccccccccc}
	\hline
	Phase   & $\chi N$ & $\eta N$ & $f$ & $\lambda$ & $\calB$ & $ \beta
	$&$N_{\br}$&$N_{\theta}$&$N_{\varphi}$&$N_{s}$\\
	\hline
	Nematic    & 6  & 12 & 0.6 & 100 & $ 1.31 $ & 2 & 32 & 8 & 17 & 100\\
	Smectic-C  & 16 & 8  & 0.6 & 100 & $ 1.31 $ & 2 & 32 & 8 & 17 & 100\\
	\hline
\end{tabular}
\label{1D:sn}
\end{table}
Firstly, the model parameters used for the self-assembly of these two structures are listed in Tab.\,\ref{1D:sn}. Ordered phases obtained through our calculation are shown in  Fig.\,\ref{Fig:nematic}. The density distributions $\phi_{A}=f,\,\phi_{B}=1-f$, and the orientation distributions $\phi_B(\bu, z=0)$ presents parallel alignment of semiflexible blocks. 
These characteristics satisfy the unique properties of the nematic phase, and its molecular schematic is shown in Fig.\,\ref{Fig:nematic} (1c).
The density distribution and the orientation distribution peaked at $\theta \in (0, \pi/2)$ in Fig.\,\ref{Fig:nematic} (2a,2b) indicate that the ordered structure corresponds to the smectic-C phase. The molecular schematic of smectic-C is presented in Fig.\,\ref{Fig:nematic} (2c).

\begin{figure}[H]
\centering
\includegraphics[width=12cm]{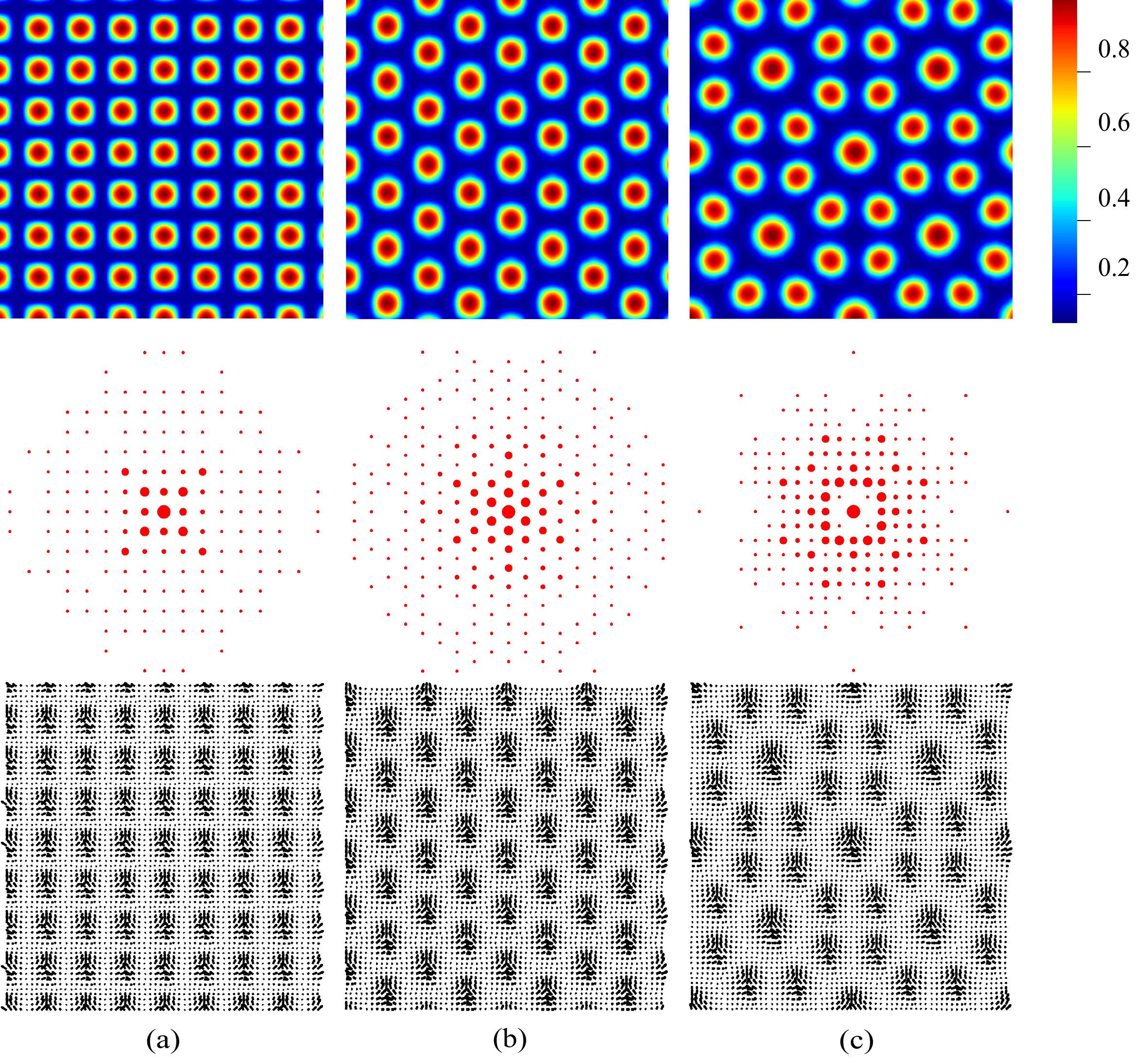}
\label{2Ddraw:ufield1}
\caption{Density distributions of monomer $B$ (first row),  FFT patterns (second
	row), and orientations distributions  of semiflexible chains (third row). (a) $[4^4]$, (b) $[3^6]$, (c) $[3^4.4;3^8]$ tilings.}
\label{2Ddraw1}
\end{figure}

In 5D calculations, flexible-semiflexible diblock copolymers can self-assemble into various morphologies of columnar structures. 
Here we investigate $[4^4]$ and $[3^6]$ tilings with the Archimedean tiling naming rules. 
In our calculations, monomers $B$ at higher concentrations form vertices, which are then connected to form polygonal tilings. Choosing a vertex as the center, we list the number of edges of each polygon connected to this vertex in a clockwise or counterclockwise order. For instance, $[4^4]$ tiling denotes a cyclic sequence of four squares. The "$[]"$ are used to distinguish these tiling pattern naming from regular numbers.
Here we not only reproduce $[4^4]$ and $[3^6]$ tilings but also discover a more complex $[3^4.4;3^8]$ tiling, with model parameters detailed in Tab.\,\ref{2D:parameter}. Furthermore, density distribution of monomer $B$, FFT patterns and orientation distributions are shown in Fig.\,\ref{2Ddraw1}.

\begin{table}[H]
\centering
\caption{The parameters of $[4^4]$, $[3^6]$ and $[3^4.4;3^8]$ tilings.}
\begin{tabular}{ccccccccccc}
	\hline
	Phase & $\chi N$ & $\eta N$ & $f$ &$\lambda$&$\calB$&$ \beta
	$&$N_{\br}$&$N_{\theta}$&$N_{\varphi}$&$N_{s}$\\
	\hline
	$[4^4]$ &20&10 & 0.75&300&$[0,1.04)^2$&2&20&4&9&100\\
	$[3^6]$ & 18&9& 0.75&300&$[0,0.52)^2$&2&20&4&9&100\\
	$[3^4.4;3^8]$ &19 & 9.5 & 0.75&300&$[0,0.78)^2$&2&20&4&9&100\\
	\hline
\end{tabular}
\label{2D:parameter}
\end{table}

\begin{figure}[H]
\centering
\includegraphics[width=14cm]{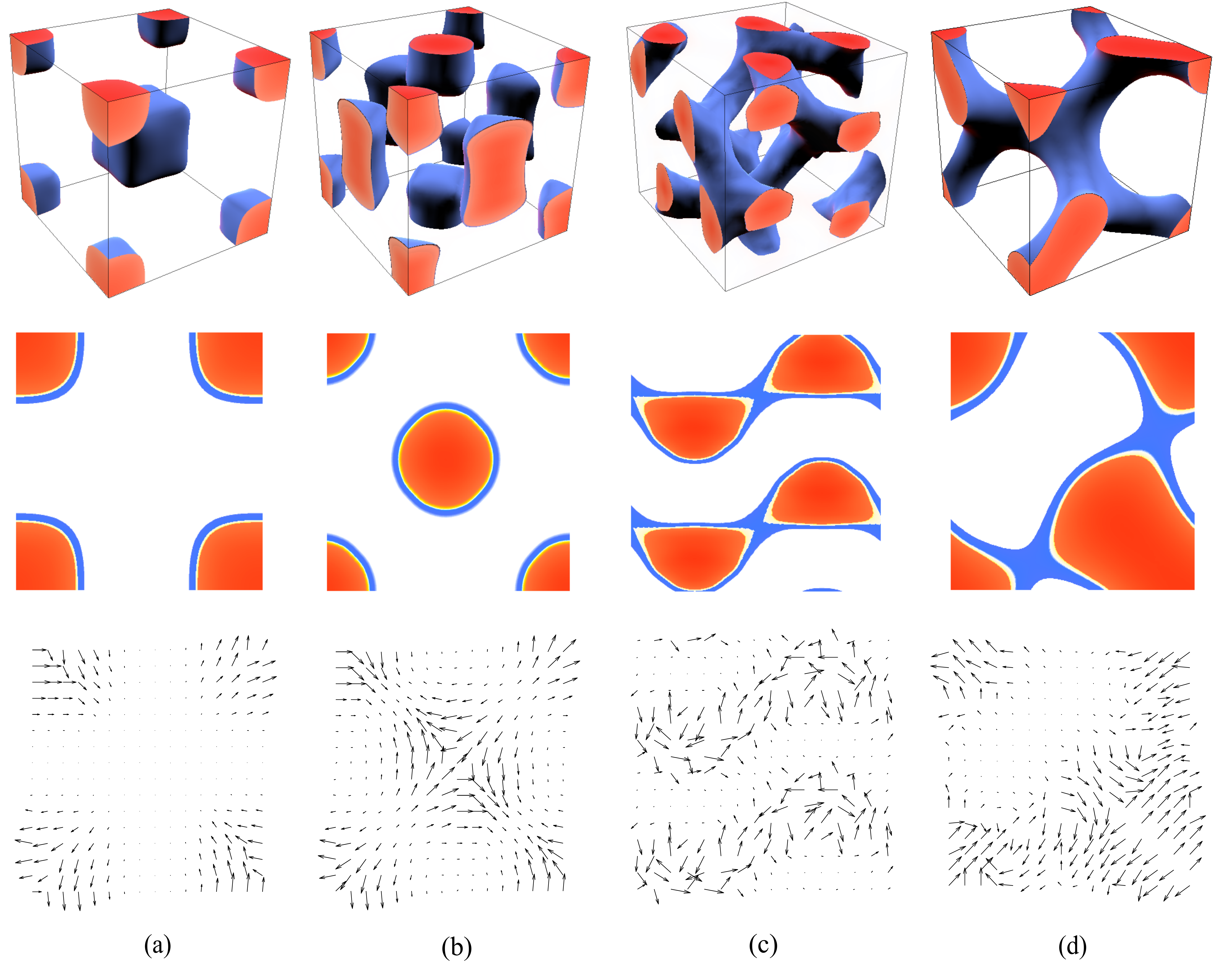}
\caption{The morphologies of phases (the first row), the top views for morphologies (second row), and the corresponding orientation distributions (the third row) of (a) BCC, (b) FCC, (c) DG and (d) SG, respectively.}
\label{3Ddraw}
\end{figure}

3D space cubic structures have been successfully predicted in rigid polymers\,\cite{ khandpur1995polyisoprene, jiang2013discovery}. Based on this, we extend our calculations to 6D case. By selecting appropriate parameters (Tab.\,\ref{3D:parameter}), we successfully simulate body-centred cubic (BCC), face-centred cubic (FCC), double gyroid (DG), and single gyroid (SG) phases.
Their morphologies in 3D space, top views of the density distributions and orientation distributions of theses phases are presented in Fig.\,\ref{3Ddraw}.

\begin{table}[H]
\centering
\caption{The parameters for ordered structures.}
\begin{tabular}{ccccccccccc}
	\hline
	Phase & $\chi N$ & $\eta N$ & $f$ &$\lambda$&$\calB$&$ \beta
	$&$N_{\br}$&$N_{\theta}$&$N_{\varphi}$&$N_{s}$\\
	\hline
	DG  & 18 & 9  &  0.6 & 100 & $[0,1.31)^3 $ & 2 & 16 & 8 & 17 & 100\\
	FCC & 14 & 14 &  0.7 & 100 & $[0,1.31)^3$ & 2 & 16 & 8 & 17 & 100\\
	BCC & 16 & 8  &  0.7 & 100 & $[0,1.31)^3 $ & 2 & 16 & 8 & 17 & 100\\
	SG  & 14 & 14 &  0.6 & 100 & $[0,1.31)^3 $ & 2 & 16 & 8 & 17 & 100\\
	\hline
\end{tabular}
\label{3D:parameter}
\end{table}

\subsection{The efficiency of optimizing the computational domain}

In this section, we express the efficiency of optimizing computational domain, which can automatically adjust
the size and shape of the computational box in the SCFT calculation. The error tolerance $\varepsilon_H$ is set as $10^{-4}$.

To show that the computational box can be adjusted automatically, we present two cases with different initial values. Fix model parameters $\chi N=18$, $f=0.75$, $\eta N=9$, $\lambda=300$, $\beta=2$ and  discretization nodes $N_{\br}=16$, $N_{\theta}=8$, $N_{\varphi}=0$, $N_{s}=100$, We illustrate the effectiveness of optimizing computational domain from the aspects of the angle $\theta$ between
edges, the length $a_1$, $a_2$ of the computational domain and the free energy H during the iteration of the optimal computational domain. 

\begin{figure}[H]
\centering
{
	\includegraphics[width=16cm]{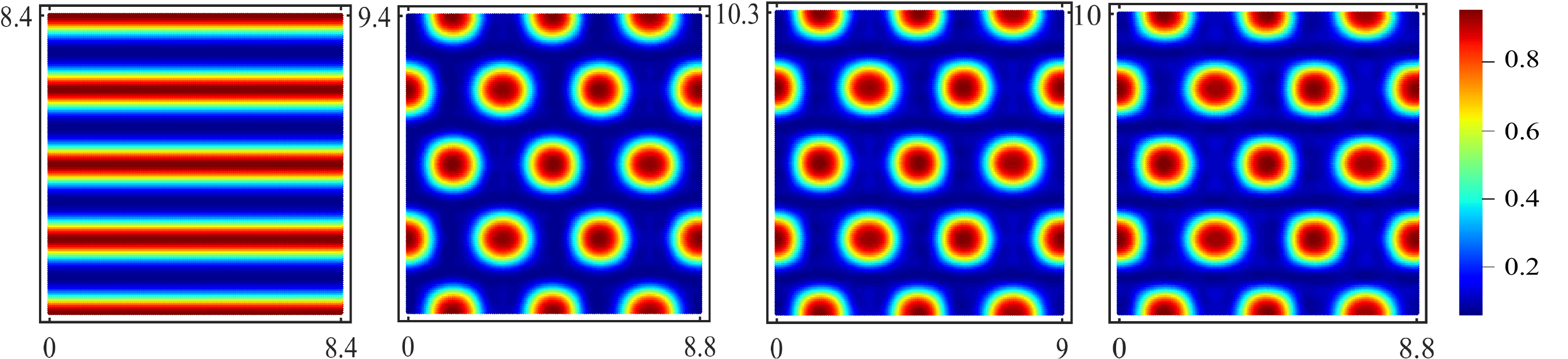}
}
\caption{The process of SCFT iteration and optimization starting from the lamellar phase as initial values, transitioning into a $[3^6]$  phase, is shown from left to right.
	The model parameters $\chi N=18$, $f=0.75$, $\eta N=9$, $\lambda=300$,
	$\beta=2$, discretization nodes $N_{\br}=20$, $N_{\theta}=8$,
	$N_{\varphi}=0$, $N_{s}=100$.}
\label{box:phase4}
\end{figure}

\begin{figure}[H]
\centering
{
	\includegraphics[width=13cm]{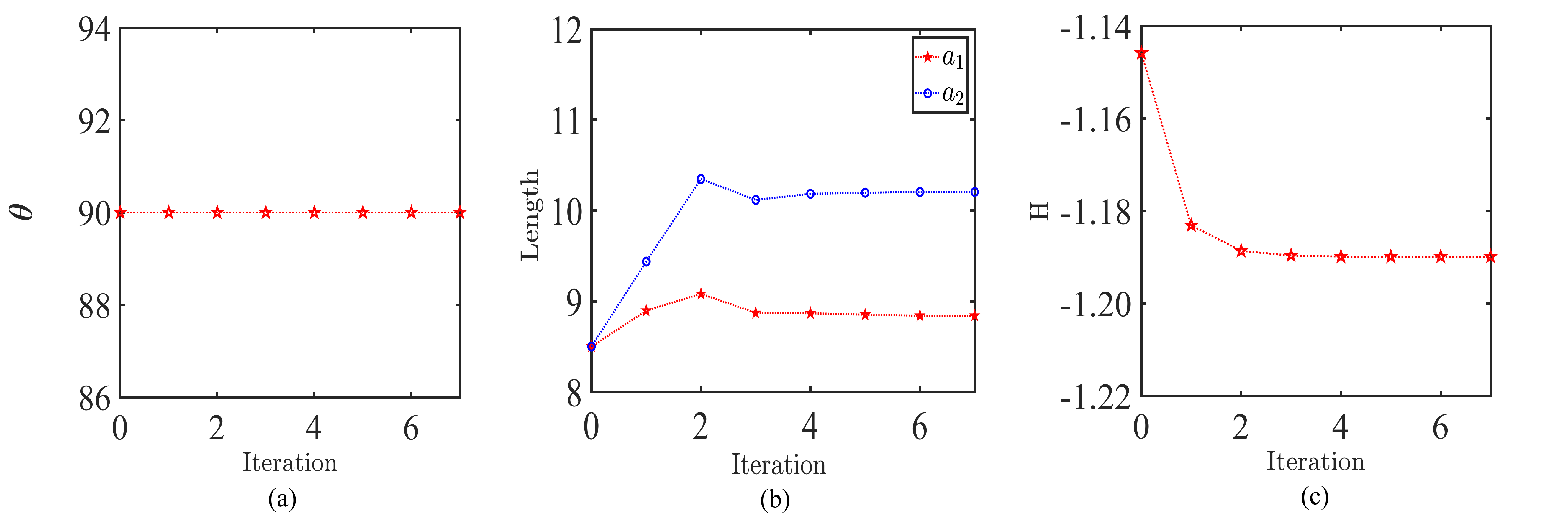}
}
\caption{The change process of (a) the angle $\theta$, (b) the box lengths $a_1$ and $a_2$, and (c) the free energy $H$ when lamellar transitions to $[3^6]$.}
\label{box:iteration3}
\end{figure}

\begin{figure}[H]
\centering
{
	\includegraphics[width=14cm]{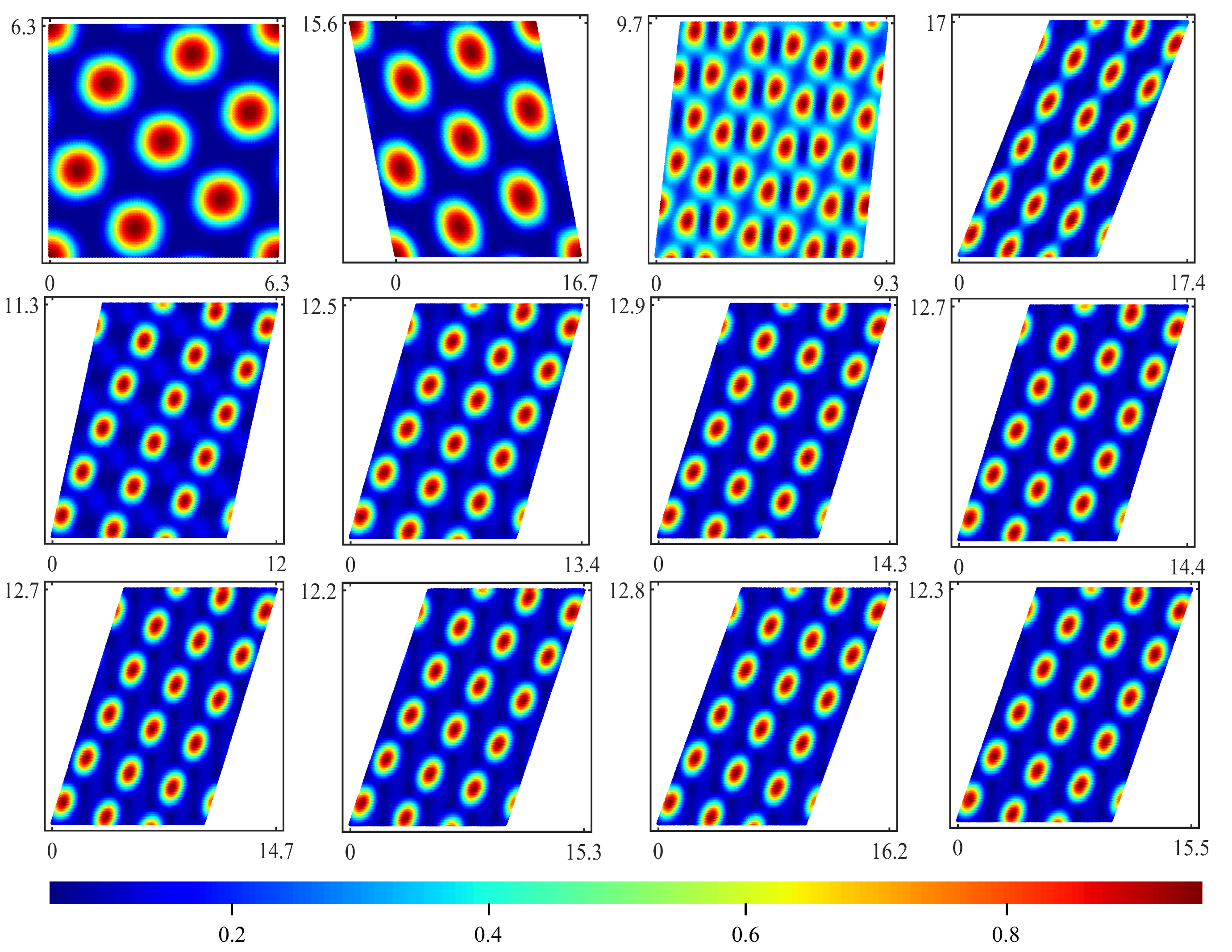}
}
\caption{The process of SCFT iteration and optimization starting from the $[3^6]$ phase as  initial values, transitioning into a $[4^4]$  phase, is shown from left to right.
	The model parameters $\chi N=18$, $f=0.75$, $\eta N=9$, $\lambda=300$, $\beta=2$; discretization nodes $N_{\br}=16$, $N_{\theta}=8$, $N_{\varphi}=0$, $N_{s}=100$.}
\label{box:phase26}
\end{figure}

\begin{figure}[H]
\centering
{
	\includegraphics[width=13cm]{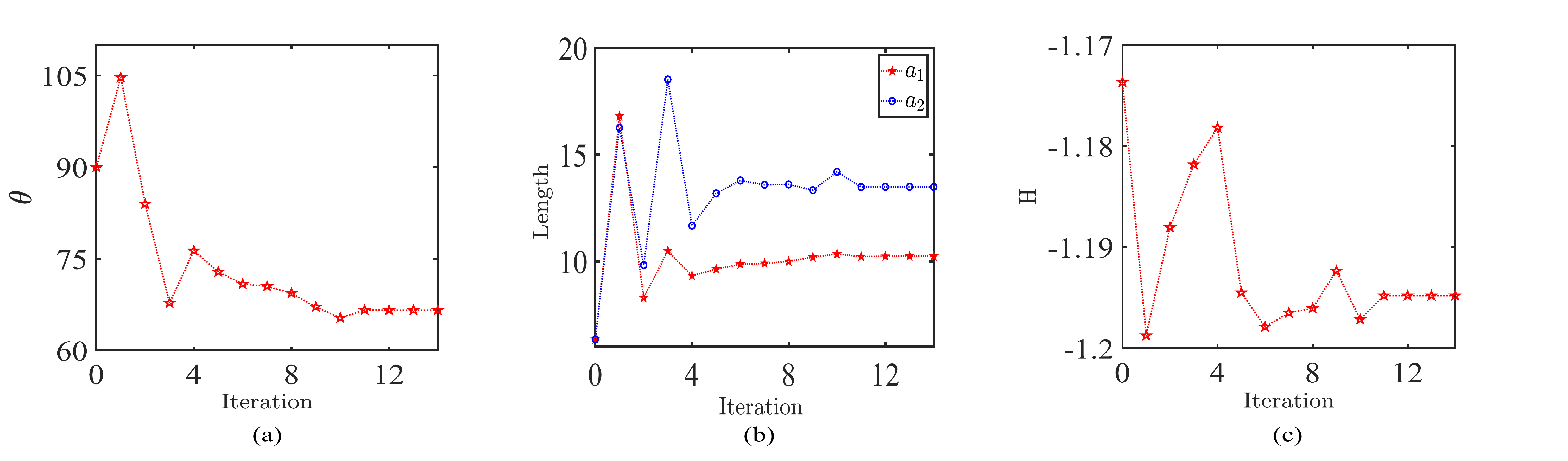}
}
\caption{The change process of (a) the angle $\theta$, (b) the box lengths $a_1$ and $a_2$, and (c) the free energy $H$ when $[3^6]$ transitions to $[4^4]$ phases.}
\label{box:iteration5}
\end{figure}

Using the lamellar phase within a square domain as the initial values, through the SCFT iteration and optimization process, the lamellar phase transitions into a $[3^6]$ phase, as shown in Fig.\,\ref{box:phase4}. The variations in $\theta$, $a_1$, $a_2$, and $H$ during the optimization process are presented in Fig.\,\ref{box:iteration3}. The constant angles and changing side lengths indicate that the optimal domain for the $[3^6]$ phase should be rectangular rather than square. The convergence of the free energy ensures the stability of the entire optimization process.
Starting from the  $[3^6]$ phase with rectangular domain, the phase transitions into a  $[4^4]$ phase with parallelogram domain, as shown in Fig.\,\ref{box:phase26}. This result is consistent with previous findings\,\cite{jiang2013self}. As shown in Fig.\,\ref{box:iteration5}, the variations in angle, side lengths, and free energy provide  support for this optimization process.

\section{Conclusion and future work}
\label{sec:concl}

In this paper, we design a set of high-accurate and efficient numerical algorithms for SCFT of liquid-crystalline polymers. We present ten numerical discretization methods for solving semi-flexible propagators and compare their accuracy and computational complexity. Meanwhile, we apply the adaptive Anderson mixing method and cascadic multi-level method to dramatically improve the
convergence rate of the Anderson mixing scheme, resulting in a notable reduction of CPU time by $63.88\%$.
Moreover, we optimize the computational domain during the SCFT calculations, successfully identifying the optimal periods for different phases and achieving more accurate free energy for these phases.
Furthermore, we investigate the self-assembly structures of flexible-semiflexible diblock copolymers in 4D, 5D, and 6D spaces.  Numerical methods proposed in this paper are adaptable to all liquid crystalline polymeric systems. In the future, we aim to apply these numerical methods to more complex liquid crystalline polymers.



\section*{Appendix} 
\begin{appendix}
\section{Spherical harmonic expansion}
\label{subsec:spe}
The square integrable function $q(\br,\bu,s)$ on unit sphere $S$($\bu\in S$) can be expanded by spherical harmonic series as
\begin{align}
	q(\br,\bu,s)=\sum_{l=0}^{\infty}\sum_{m=-l}^{m=l}q_{lm}(\br,s)Y_{l,m}(\bu),
	\label{sphrical:harmomic}
\end{align}
with $Y_{l,m}(\bu)$ is the spherical harmonic
function\,\cite{fredrickson2006equilibrium}, and
\begin{align}
	\begin{split}
		q_{lm}(\br,s)=&\int_S q(\br,\bu,s)Y^{*}_{l,m}(\bu)~d\bu\\
		=&\int_{0}^{\pi}\int_{0}^{2\pi}\sin\theta
		q(\br,\theta,\varphi,s)Y^{*}_{l,m}(\theta,\varphi)~d\varphi \,d\theta,
		\label{harmonic:cofficient}
	\end{split}
\end{align}
where $Y^*_{l,m}$ is the complex conjugate of $Y_{l,m}$.
The vector $\bu$ can be expressed through unit sphere coordinates
\begin{align}
	\bu=(x,y,z)^{T}=(\cos\varphi\sin\theta,\sin\theta\sin\varphi,\cos\theta)^{T},
	~~~~~~~\varphi\in [0, 2\pi),~\theta\in [0,\pi).
\end{align}
The spherical harmonic function $Y_{l,m}(\bu)$ can be expressed as
\begin{align}
	Y_{l,m}(\bu)=Y_{l,m}(\theta,\varphi)=(-1)^m\left[\frac{2l+1}{4\pi}\frac{(l-m)!}{(l+m)!}
	\right]^{1/2} P_l^m(\cos\theta)e^{im\varphi},
\end{align}
where $P_l^m(x)$ is the associated Legendre function
\begin{equation}
	\begin{aligned}
		P_{l}^{m}(x)&=(1-x^2)^{\frac{m}{2}}\frac{d^m}{dx^m}P_l\left(x\right),\quad l,\;m\geqslant 0,\\
		P_{l}^{m}\left(x\right)&=\left(-1\right)^m\frac{\left(l-\left| m\right|\right)!}{\left(l+\left| m\right|\right)!}P_{l}^{\left| m\right|}\left(x\right),\quad m<0,
	\end{aligned}
\end{equation}
with $P_l(x)$ is the Legendre polynomial
\begin{align}
	P_l(x)=\frac{1}{2^ll!}\frac{d^l}{dx^l}{(x^2-1)^l},{\qquad}l\ge0.
\end{align}
Besides, the spherical harmonic function $Y_{l,m}(\bu)$ is the
characteristic function of the Laplace operator $\nabla^{2}_{\bu}$  on the unit
sphere,
\begin{align}
	\nabla^{2}_{\bu} Y_{l,m}(\bu)=-l(l+1)Y_{l,m}(\bu).
\end{align}
The spherical harmonic function is a complete orthogonal basis on the sphere, which satisfies orthogonality
\begin{equation}
	\begin{aligned}
		\int_S \left[Y_{l,m}(\bu)\right]^{*} Y_{l^{\prime},m^{\prime}}(\bu) d \bu
		&=\int_{0}^{2 \pi} \int_{0}^{\pi} \sin \theta\left[Y_{l,m}(\theta,
		\varphi)\right]^{*} Y_{l^{\prime},m^{\prime}}(\theta, \varphi) d \theta d
		\varphi  \\
		&=\delta_{l l^{\prime}} \delta_{m m^{\prime}}.
	\end{aligned}
\end{equation}
The SHT and ISHT can be realized by the SPHEREPACK function package\,\cite{Adams1999SPHEREPACK}.

\section{Fourier pseudo-spectral method}
\label{subsec:FFT}
The Fourier pseudo-spectral method is a class of numerical methods for solving partial differential
equations used in
applied mathematics and scientific computing. It is closely related to the Fourier spectral method but complements it with an additional pseudo-spectral basis, representing functions on a
quadrature grid. The Fourier pseudo-spectral method simplifies the evaluation
of certain operators and can considerably speed up the calculation with FFT. The Fourier pseudo-spectral method is a natural
choice for obtaining optimal spatial accuracy under periodic boundary
conditions.

Before we show the details of the Fourier pseudo-spectral methods, a brief introduction to the Bravais lattice and reciprocal lattice is necessary.
The Bravais lattice is defined by $R_n =n_1\ba_1 +\cdots +n_d\ba_d$, $n_1,\cdots,n_d \in \mathbb{Z}$, the
primitive vectors are denoted by $\ba_i = (a_{i1}, \cdots, a_{id})$, where $ i = 1,~\cdots ,~d $. $d$ is the
geometric dimension of space. Reciprocal vectors $\calB =
(\bb_1,\cdots,\bb_d)$. The corresponding reciprocal lattice primitive
vectors are $\bb_i =(b_{i1}, \cdots, b_{id})$, $\Lambda^*=k_1\bb_1+\cdots
+k_{d}\bb_{d},k_1,\cdots,k_d \in \mathbb{Z}$ is the reciprocal lattice. These two sets of primitive vectors satisfy $\ba_i \cdot \bb_j
= 2\pi \delta_{ij}$, where $ i,~j = 1,~\cdots ,~d $. 

Taking the 1D spatial case as an example, the spatial variable $r$ is discretized into grid points $r_j$, $j=1,\cdots, N$. The periodic function can be expanded as\,\cite{Jiang2010spectral}
\begin{align}
	f(r_j)=\sum_{k=-N}^N\hat{f}(k) e^{i k\cdot r_j},
	\label{Fourier:expansion}
\end{align}
where the discrete Fourier coefficient $\hat{f}$ can be obtained by FFT. Applying the Fourier pseudo-spectral method to solve the propagator (Eqn.\,\eqref{simple:PDE}), it can be expanded into
\begin{align}                      
	\frac{\partial}{\partial s}\sum_{k = -N}^{N}\hat{q}_{A}(k,s)e^{ik\cdot r}
	=\nabla^{2}_{r}\sum_{k =
		-N}^{N}\hat{q}_{A}(k,s)e^{i k\cdot r}-\sum_{k_{1}+k_{2}=-N}^{N}\hat{w}_A(k_{1})e^{i k_{1}\cdot r}\sum_{k_{2}=
		-N}^{N}\hat{q}_{A}(k_{2},s)e^{i k_{2}\cdot r},
	\label{FFT:expan1}                                                          
\end{align}                                                                     
Take the inner product of both sides of Eqn.\,\eqref{FFT:expan1} with
$e^{i k\cdot r}$, we get
\begin{align}                                                                   
	\frac{\partial}{\partial
		s}\hat{q}_{A}(k,s)=-k^2\hat{q}_{A}(k,s)-\sum_{k_{1}+k_{2}=-N}^{N}\hat{w}_A(k_{1})\hat{q}_{A}(k_{2},s).
	\label{FFT:expan2}                                                          
\end{align}
Since the last term of the above formula creates a convolution that makes the
calculation more complicated, we set $G(r,s) = w(r)q(r,s)$,
and get the semi-discrete scheme
\begin{align}                                                                   
	\frac{\partial}{\partial
		s}\hat{q}_{A}(k,s)=-k^2\hat{q}_{A}(k,s)-\hat{G}(k,s).
	\label{FFT:expan3}
\end{align}                       
FFT can be performed for the current calculation.
Since propagators $q_{A}(\br,s)$, $q^{\dagger}_{A}(\br,s)$,
$q_{B}(\br,\bu,s)$ and $q^{\dagger}_{B}(\br,\bu,s)$ (see Eqns.\,\eqref{propagator:qA:PDE}-\eqref{inverse:qA:PDE}) satisfy the periodic boundary conditions  on the spatial variable, we use the Fourier pseudo-spectral method to solve them.

\section{Full discretization scheme of the four-dimensional semiflexible propagator equation}
\label{app:full}

Here we present the full discretization scheme with 1-dimensional space variable $r$, 2-dimensional orientation variable $\bu =(\theta,\varphi)$, and 1-dimensional time variable $s$ for equation
\begin{equation}
	\frac{\partial}{\partial s} q_{B}(r,\bu,s)=-\beta \cos\theta\, \nabla_{r}
	q_{B}(r,\bu,s)+\Gamma(r,\bu)q_{B}(r,\bu,s)+\frac{1}{2\lambda}\nabla^2_\bu q_{B}(r,\bu,s),~~~~~ 1-f\leq s\leq 1,
	\label{fulusan}
\end{equation}
where $\Gamma(r,\bu)=-\big(w_{B}(r)-\bM(r):\bigr[\bu\bu-\frac{1}{3}\bI\bigm]\bigl)$.
By expanding $q_B(r,\bu,s)$ in a Fourier series with respect to the spatial variable $r$ and in a spherical harmonic series with respect to the orientation variable $\bu$, the equation can be formulated as
\begin{equation}
	\sum_{k=-N}^{N-1} \sum_{l=0}^{L} \sum_{m=-l}^{l}\frac{\partial
		\tilde{q}_B(k,l,m,s)}{\partial s}e^{ik r}Y_{l,m}(\bu)=M_1+M_2+M_3,
	\label{ODE2}
\end{equation}
where $2N+1$ is the number of discrete points in space, $L + 1$ is the number of
discrete nodes of $\theta$, and 
\begin{equation}
	\begin{aligned}
		M_1 &=\sum_{k=-N}^{N-1} \sum_{l=0}^{L} \sum_{m=-l}^{l}
		-\frac{l(l+1)}{2\lambda}\tilde{q}_B(k,l,m,s)e^{ikr}Y_{l,m}(\bu),\\
		M_2 &=\sum_{k=-N}^{N-1}\sum_{l=0}^{L}\sum_{m=-l}^{l}\tilde{\Gamma}_B(k,l,m,s)e^{ik
			r}Y_{l,m}(\bu),\\
		M_3 &=\sum_{k=-N}^{N-1}\sum_{l=0}^{L}\sum_{m=-l}^{l}-i \beta
		\tilde{U}_B(k,l,m,s)e^{ikr}Y_{l,m}(\bu).
	\end{aligned}
\end{equation}
The coefficients 	$\tilde{q}_B(k,l,m,s)$, $\tilde{U}_B(k,l,m,s)$, $\tilde{\Gamma}_B(k,l,m,s)$ are obtained by first performing Fourier transforms on $ q_{B}(r,\bu,s)$,
$q_{B}(r,\bu,s)kcos\theta$ and $\Gamma(r,\bu)q_{B}(r,\bu,s)$, followed by spherical harmonic transforms, respectively.
Using the testing function $v=e^{ik'r}Y_{l',m'}(\theta,\varphi)$, Eqn.\eqref{ODE2} can be written as

\begin{equation}
	\begin{split}
		&\sum_{p=1}^{2N}\sum_{k=-N}^{N-1} \sum_{j=1}^{L+1}\sum_{l=0}^{L} \sum_{n=1}^{2L+1}\sum_{m=-l}^{l}
		\frac{\partial \tilde{q}_B(k,l,m,s)}{\partial s}
		e^{ikr_{p}}e^{-ik'r_{p}}Y_{l,m}(\theta_{j},\varphi_{n})Y^{*}_{l',m'}(\theta_{j},\varphi_{n})\\
		&=\sum_{p=1}^{2N}\sum_{k=-N}^{N-1} \sum_{j=1}^{L+1}\sum_{l=0}^{L} \sum_{n=1}^{2L+1}\sum_{m=-l}^{l} -\frac{l(l+1)}{2\lambda}
		\tilde{q}_B(k,l,m,s)e^{ikr_{p}}e^{-ik'r_{p}}Y_{l,m}(\theta_{j},\varphi_{n})Y^{*}_{l',m'}(\theta_{j},\varphi_{n})\\
		&+\sum_{p=1}^{2N}\sum_{k=-N}^{N-1} \sum_{j=1}^{L+1}\sum_{l=0}^{L}
		\sum_{n=1}^{2L+1}\sum_{m=-l}^{l}\tilde{\Gamma}_B(k,l,m,s)e^{ik
			r_{p}}e^{-ik'r_{p}}Y_{l,m}(\theta_{j},\varphi_{n})Y^{*}_{l',m'}(\theta_{j},\varphi_{n})\\
		&-\sum_{p=1}^{2N}\sum_{k=-N}^{N-1} \sum_{j=1}^{L+1}\sum_{l=0}^{L} \sum_{n=1}^{2L+1}\sum_{m=-l}^{l}
		i \beta  \tilde{U}_B(k,l,m,s)e^{ik
			r_{p}}e^{-ik'r_{p}}Y_{l,m}(\theta_{j},\varphi_{n})Y^{*}_{l',m'}(\theta_j,\varphi_{n}),
		\label{ODE5}
	\end{split}
\end{equation}
where $\theta_{j}=\frac{\pi j}{L+1}, j=1,2\dots L+1$, $\varphi_{n}=\frac{2\pi
	n}{2L+1}, n = 1,2\dots 2L+1$, $r_{p}=\frac{|\calB|p}{2N}, p=1,2\dots 2N$, and the periodic boundary condition makes $q_B(r_0,\bu,s) = q_B(r_{2N}, \bu,s)$.

Define the discrete form of inner-product as
\begin{align}
	(q_B,v)=\sum_{p=1}^{2N}\sum_{j=1}^{L+1}\sum_{n=1}^{2L+1}q_B(r_p,\theta_{j},\varphi_{n})v^*(r_p,\theta_j,\varphi_{n}),
\end{align}
where $v^*$ is the complex conjugate of $v$. Thus, Eqn.\eqref{ODE5} can be written as
\begin{equation}
	\begin{split}
		&\sum_{k=-N}^{N-1} \sum_{l=0}^{L} \sum_{m=-l}^{l}\frac{\partial
			\tilde{q}_B(k,l,m,s)}{\partial s}(e^{ik r}Y_{l,m}(\theta,\varphi),e^{ik'
			r}Y_{l',m'}(\theta,\varphi))\\
		&=\sum_{k=-N}^{N-1} \sum_{l=0}^{L} \sum_{m=-l}^{l}
		-\frac{l(l+1)}{2\lambda}\tilde{q}_B(k,l,m,s)(e^{ik
			r}Y_{l,m}(\theta,\varphi),e^{ik' r}Y_{l',m'}(\theta,\varphi))\\
		&+\sum_{k=-N}^{N-1}\sum_{l=0}^{L}\sum_{m=-l}^{l}\tilde{\Gamma}_B(k,l,m,s)(e^{ik
			r}Y_{l,m}(\theta,\varphi),e^{ik' r}Y_{l',m'}(\theta,\varphi))\\
		&-\sum_{k=-N}^{N-1}\sum_{l=0}^{L}\sum_{m=-l}^{l}i \beta
		\tilde{U}_B(k,l,m,s)(e^{ik r}Y_{l,m}(\theta,\varphi),e^{ik'
			r}Y_{l',m'}(\theta,\varphi)).
		\label{ODE3}
	\end{split}
\end{equation}
Simplifying it leads to the full discretization form 
\begin{equation}
	\begin{split}
		\frac{\partial \tilde{q}_B(k,l,m,s)}{\partial s}
		=-\frac{l(l+1)}{2\lambda}\tilde{q}_B(k,l,m,s)+\tilde{\Gamma}_B(k,l,m,s)-i \beta
		\tilde{U}_B(k,l,m,s).
		\label{ODE4}
	\end{split}
\end{equation}
\end{appendix}

\bibliographystyle{unsrt} 
\bibliography{ref1} 

\end{document}